\documentclass[journal]{IEEEtran}

\usepackage{graphicx}

\usepackage[labelformat=simple]{subcaption}

\usepackage{xcolor}
\usepackage{amsmath} 
\usepackage{mathrsfs}
\usepackage{amssymb}
\usepackage{mathtools}
\usepackage{mathdots}
\usepackage{cite}
\usepackage[hidelinks]{hyperref}
\let\labelindent\relax
\usepackage[shortlabels]{enumitem}
\usepackage{algorithm}
\usepackage{algorithmic}
\usepackage{bm}
\usepackage{url}
\usepackage{optidef}
\usepackage{multirow}
\usepackage{dsfont}
\usepackage{stfloats}
\usepackage{caption}

\DeclareSymbolFont{bbold}{U}{bbold}{m}{n}
\DeclareSymbolFontAlphabet{\mathbbold}{bbold}

\newcommand{\vect}[1]{\mathbbold{#1}}

\newcommand{\vectorzeros}[1][]{\vect{0}_{#1}}

\usepackage[normalem]{ulem}
\newcommand\rout{\bgroup\markoverwith{\textcolor{red}{/}}\ULon} 

\usepackage[prependcaption,colorinlistoftodos]{todonotes}

\newtheorem{theorem}{\bf Theorem}[section]
\newtheorem{proposition}{\bf Proposition}[section]
\newtheorem{lemma}{\bf Lemma}[section]

\newtheorem{definition}{\bf Definition}[section]
\newtheorem{remark}{\bf Remark}[section]

\newcommand{\longthmtitle}[1]{\mbox{}\textup{\textbf{(#1):}}}

\newcommand{\real}{{\mathbb{R}}}

\newcommand{\realpositive}{{\mathbb{R}}_{>0}}
\newcommand{\realnonnegative}{{\mathbb{R}}_{\ge 0}}
\newcommand{\integers}{\mathbb{Z}}
\newcommand{\integersnonnegative}{\mathbb{Z}_{\geq 0}}
\newcommand{\integerspositive}{\mathbb{Z}_{> 0}}

\newcommand{\eps}{\epsilon}

\newcommand{\oprocendsymbol}{\hbox{$\bullet$}}
\newcommand{\oprocend}{\relax\ifmmode\else\unskip\hfill\fi\oprocendsymbol}

\newcommand{\defeq}{\vcentcolon=}

\newcommand{\col}{{\textup{col}}}
\newcommand{\Prob}{\mathbb{P}}
\newcommand{\Exp}{\mathbb{E}}

\newcommand{\bv}{{\mathscr{B}}}
\newcommand{\Tini}{{T_{\textup{ini}}}}
\newcommand{\Tf}{{T_{\textup{f}}}}
\newcommand{\uini}{{u_{\textup{ini}}}}
\newcommand{\uinihat}{{\hat{u}_{\textup{ini}}}}
\newcommand{\yini}{{y_{\textup{ini}}}}
\newcommand{\yinihat}{{\hat{y}_{\textup{ini}}}}

\newcommand{\xini}{{x_{\textup{ini}}}}
\newcommand{\xinihat}{{\hat{x}_{\textup{ini}}}}

\newcommand{\Yp}{{Y_{\mathrm{p}}}}
\newcommand{\Yf}{{Y_{\mathrm{f}}}}
\newcommand{\Uphat}{{\widehat{U}_{\mathrm{p}}}}
\newcommand{\Ufhat}{{\widehat{U}_{\mathrm{f}}}}
\newcommand{\Yphat}{{\widehat{Y}_{\mathrm{p}}}}
\newcommand{\Yfhat}{{\widehat{Y}_{\mathrm{f}}}}

\newcommand{\Phat}{{\widehat{\Prob}}}
\newcommand{\B}{{B_{\eps}(\Phat_N)}}
\newcommand{\Q}{\mathbb{Q}}
\newcommand{\ambiguity}{\widehat{\mathcal{P}}}

\begin{document}

\title{Distributionally Robust Chance Constrained Data-enabled Predictive Control}
\author{Jeremy Coulson \qquad John Lygeros \qquad Florian D\"{o}rfler%
\\\medskip\vspace{0.08cm}
\fbox{
\small{Accepted for publication: IEEE Transactions on Automatic Control. (DOI: 10.1109/TAC.2021.3097706) \textcopyright~2021 IEEE.}
}%
\vspace{-1.0cm}%
\thanks{All authors are with the Department of Information Technology and Electrical Engineering at ETH Z\"{u}rich, Switzerland~\texttt{\{jcoulson, lygeros, dorfler\}@control.ee.ethz.ch}. The research was supported by the ERC under project OCAL, grant agreement 787845, and ETH Z\"urich funds.}
}


%

\maketitle

\begin{abstract}
We study the problem of finite-time constrained optimal control of 
unknown stochastic linear time-invariant systems, which is the key ingredient of a predictive control algorithm -- albeit typically having access to a model. We propose a novel distributionally robust data-enabled predictive control (DeePC) algorithm which uses noise-corrupted input/output data to predict future trajectories and compute optimal control inputs while satisfying output chance constraints. The algorithm is based on (i) a non-parametric representation of the subspace spanning the system behaviour, where past trajectories are sorted in Page or Hankel matrices; and (ii) a distributionally robust optimization formulation which gives rise to strong probabilistic performance guarantees. We show that for certain objective functions, DeePC exhibits strong out-of-sample performance, and at the same time respects constraints with high probability. The algorithm provides an end-to-end approach to control design for unknown stochastic linear time-invariant systems. We illustrate the closed-loop performance of the DeePC in an aerial robotics 
case study.%
\end{abstract}



\section{Introduction}
\IEEEPARstart{O}{ptimal} control of unknown systems can be approached in two ways: model-based and data-driven. In model-based control, a predictive model for the system of interest is first identified from data and subsequently used for control design. What have come to be known as ``data-driven methods'', on the other hand, aim to design controllers directly from data, without explicitly identifying a predictive model. These methods are suitable for applications where first-principle models are not conceivable (e.g., in human-in-the-loop applications), when models are too complex for control design (e.g., in fluid dynamics), and when thorough modelling and parameter identification is too costly (e.g., in robotics). 
Data-driven control has recently gained a lot of popularity, but most methods cannot be applied (respectively, lack formal certificates)
for real-time control of safety-critical systems.

In this work, we focus on a data-driven control technique for unknown, stochastic, and constrained linear systems. In particular, we present a method for finite-horizon optimal predictive control using input/output data from the unknown system, where the system behaviour is characterized by a data matrix time series.
This method was first presented for deterministic systems in~\cite{JC-JL-FD:18} and was later extended to stochastic systems in~\cite{JC-JL-FD:19}. These works were motivated by~\cite{IM-PR:08} in which a unique method of direct data-driven open-loop control was conceived based on the seminal work on behavioural systems theory~\cite{JCW-PR-IM-BDM:05}.

A key challenge for such data-driven control methods is ensuring the performance and safety of the system in the presence of uncertainties, corrupted data, and noise. In this work, we present an end-to-end optimal data-driven control approach which comes with such guarantees but without explicitly identifying a predictive model and is agnostic to the particular probabilistic uncertainty. The approach combines the so-called data-enabled predictive control (DeePC) method from~\cite{JC-JL-FD:18,JC-JL-FD:19} and distributionally robust optimization techniques from~\cite{PME-DK:18}.

Data-driven control has been historically approached using, e.g., iterative feedback tuning and virtual reference feedback tuning~\cite{HH-MG-SG-OL:98, MC-AL-SS:02}. More recently the default approach is often reinforcement learning. There are many approaches to reinforcement learning, many of which are episodic, where learning and control alternate; see~\cite{BR:19} for a review of methods and challenges. Here, we follow different lines of literature.

\emph{Behavioural Framework for Control:} 
Instead of learning a parametric system representation, one can describe the entire subspace of possible trajectories of a linear time-invariant (LTI) system only in terms of raw data sorted into a Hankel matrix. This result became known as the fundamental lemma~\cite{JCW-PR-IM-BDM:05}, has been inspired by subspace system identification~\cite{PVO-BDM:12-book}, and was first leveraged in~\cite{IM-PR:08} for the computation of open-loop control and for simulating system responses. The result was used in~\cite{JC-JL-FD:18} in which a predictive control algorithm was proposed and was extended for stochastic systems in~\cite{JC-JL-FD:19}.
Robust closed-loop stability guarantees have been provided in~\cite{JB-JK-MM-FA:20-journal}.
Additionally, numerical case studies have illustrated that the algorithm performs robustly on some stochastic and weakly nonlinear systems and often outperforms system identification followed by conventional model predictive control (MPC)~\cite{LH-JC-JL-FD:19,LH-JC-JL-FD:20,EE-JC-PB-JL-FD:19}. The behavioural framework was used in~\cite{CDP-PT:19} to construct explicit feedback controllers. Since then, the behavioural framework has become popular for control design giving rise to numerous methods~\cite{JB-AR-CWS-FA:19,JB-JK-MM-FA:20,AK-JB-JK-FA:20,HJVW-CDP-MKC-PT:20,HVW-JE-HT-MC:20,CDP-PT:20}.

\emph{Learning-based MPC:}  
In learning-based MPC the unknown system dynamics are substituted with a learned model mapping inputs to output predictions; see~\cite{LH-KW-MM-MZ:19} for a comprehensive survey. A learning MPC approach for iterative tasks is presented in~\cite{UR-FB:18}. Other approaches conceptually related to ours exploit previously measured trajectories (i.e., motion primitives or trajectory dictionaries) based on which they synthesize new trajectories~\cite{EK-NK-SB:18,JS-DR-TA-DP-GG:19,RK-DD:15}. These methods, however, require learning a dictionary of libraries that best fits a dataset, but do not take into account the control objective.

\emph{Sequential System Identification and Control:} 
System identification (ID) produces a nominal model and an uncertainty estimate allowing for robust control design. Most approaches offer asymptotic guarantees, but some provide finite sample certificates~\cite{MCC-EW:02, MV-RLK:08, ST-RB-AP-BR:17, TS-AR-MAD:19}. In this spirit, an end-to-end ID and control pipeline is given in~\cite{RB-NM-BR:18,SD-HM-NM-BR-ST:19} and arrives at a data-driven control solution with guarantees on the sample efficiency, stability, performance, and robustness. We refer also to the earlier surveys on identification for control \cite{HH:05,MG:05}.
Well-known short-comings of sequential ID and control are as follows \cite{HH:05, BAO:96,JR:93,MG:05}: The system ID step is known to be the most time consuming part of model-based control design. Furthermore, most system ID techniques seek the model that best fits the data, but do not take into account the control objective, possibly leading to unsatisfactory performance. Finally, practitioners often demand end-to-end automated solutions. Lastly, we mention approaches related to subspace identification and (in hindsight) also to the behavioural perspective: \cite{WF-BDM-MG:99} derives a linear model from a Hankel matrix to be used for predictive control (see also~\cite{BH-RK:08-book} for an overview of subspace methods for identification and predictive control). Connections between indirect data-driven control (sequential system ID and control) and direct data-driven control are discussed in~\cite{FD-JC-IM:21}.

\emph{Stochastic MPC:} To account for the stochastic uncertainty in the system dynamics, stochastic MPC approaches usually consider minimizing the expectation of the objective function, while satisfying chance constraints~\cite{AM:16}. While chance constraints may be unsuitable for applications where constraint violation has catastrophic consequences, they are often used when transient violations are allowed or when violations are associated with financial penalties. Closest to our work is~\cite{BVP-DK-PG-MM:15} in which a stochastic MPC approach is developed for the case when the probability distribution of the stochastic disturbance is unknown. The authors consider a distributionally robust approach to safe-guard against the unknown distribution assuming that an accurate system model is given.
  
The approach of this paper follows our earlier work~\cite{JC-JL-FD:18,JC-JL-FD:19}. Our contributions are as follows:
\\
\textbf{Distributionally robust data-driven control formulation:} We formulate a novel distributionally robust data-enabled predictive control problem based on behavioural systems theory and  distributionally robust optimization. \\
\textbf{Probabilistic guarantees on performance:} We show that a distributionally robust optimal control problem admits a tractable reformulation, and with high confidence, its solutions exhibit strong out-of-sample guarantees. That is, we prove that the optimal value of the tractable formulation is an upper confidence bound on the out-of-sample performance of an optimal worst-case solution. We provide sample complexity results, and  are able to leverage additional data measurements in comparison to~\cite{JC-JL-FD:19}. Furthermore, the distributional nature of the robustness means that the method is robust against a set of systems compatible with the data collected, which can include non-Gaussian
 and non-additive noise as well as (weakly) nonlinear systems.\\
\textbf{Safety through chance constraint satisfaction:} We enforce chance constraints on the outputs of the system using a distributionally robust conditional value-at-risk formulation. We provide a tractable reformulation of the chance constraints and provide sample complexity results such that they hold with high confidence and are agnostic to the underlying probabilistic uncertainty. \\
\textbf{Tight guarantees due to new data structure}: We provide a novel formulation of the so-called fundamental lemma from behavioural system theory by proving that the subspace of trajectories of an LTI system can be spanned by trajectories organized in a Page matrix. This is in contrast with all of the literature leveraging the behavioural framework for control in which a Hankel matrix is used. In particular, this provides a novel contribution relative to~\cite{JC-JL-FD:18,JC-JL-FD:19} as well as the body of literature leveraging the behavioural framework for control. Furthermore, we show that this alternative matrix structure gives tighter guarantees and results in better performance.\\
The guarantees provided are for the finite-horizon optimal control problem. Providing guarantees for a closed-loop receding horizon implementation is a subject for future work. We refer to~\cite{JB-JK-MM-FA:20-journal} for related results in this direction.

Section~\ref{sec:behavioural} reviews preliminaries on behavioural systems and presents non-parametric representations. In Section~\ref{sec:detdeepc}, we present the data-enabled predictive control algorithm for deterministic systems. In Section~\ref{sec:stodeepc} we present a distributionally robust data-enabled predictive control algorithm for stochastic systems. In Section~\ref{sec:numericalresults}, we illustrate the main results in simulation on a quadcopter and present a detailed analysis of the hyperparameters involved. We make concluding remarks in Section~\ref{sec:conclusion}. All proofs have been deferred to the Appendix.

\emph{Notation:} We denote by $\integersnonnegative$, and $\integerspositive$ the set of non-negative and positive integers respectively. Given $x,y\in\real^n$, $\langle x,y\rangle\defeq x^{\top}y$ denotes the usual inner product on $\real^n$. We denote the associated dual norm of a norm $\|\cdot\|$ on $\real^n$ by $\|x\|_{\ast}\defeq\sup_{\|y\|\leq 1}\langle x,y\rangle$. The convex conjugate of a function $f\colon X\to\real$ is denoted by $f^{\ast}(\theta)\defeq \sup_{x\in X} \langle \theta,x \rangle-f(x)$. We define the positive part of a real-valued function $f$ as $f_{+}(x)\defeq \max\{f(x),0\}$. We denote by $\delta_x$ the Dirac distribution at $x$. Given a signal $u\colon \integers \to \real^{m}$, we denote the restriction of the signal to an interval by $u_{[1,T]}\defeq \col(u_1,\dots, u_T)$, where $\col(u_1,\dots,u_T)$ denotes the stacked column vector $(u_1^\top,\dots,u_T^\top)^\top$. We use the $\widehat{\cdot}$ symbol to denote recorded data samples and to indicate that objects depend on data samples.
\section{Behavioural Systems}\label{sec:behavioural}
\subsection{Preliminaries and Notation}
Behavioural system theory is a natural way of viewing a dynamical system when one is not concerned with a particular system representation.
This is in contrast to classical system theory, where a particular parametric system representation (such as a state-space model) is used to describe the behaviour, and system properties are derived by studying the chosen representation.
Following~\cite{IM-JCW-SVH-BDM:06-book}, we define a dynamical system and its properties in terms of its behaviour.
\begin{definition}
A \emph{dynamical system} is a $3$-tuple $(\integersnonnegative,\mathbb{W},\bv)$ where $\integersnonnegative$ is the discrete-time axis, $\mathbb{W}$ is a signal space, and $\bv\subseteq \mathbb{W}^{\integersnonnegative}$ is the behaviour.
\end{definition}
\begin{definition}
Let $(\integersnonnegative,\mathbb{W},\bv)$ be a dynamical system.
\begin{enumerate}[(i)]
\item $(\integersnonnegative,\mathbb{W},\bv)$ is \emph{linear} if $\mathbb{W}$ is a vector space and $\bv$ is a linear subspace of $\mathbb{W}^{\integersnonnegative}$.\label{item:systemdef1}
\item $(\integersnonnegative,\mathbb{W},\bv)$ is \emph{time invariant} if $\bv \subseteq \sigma\bv$ where $\sigma$ is the backward time shift defined by $(\sigma w)(t)=w(t+1)$ and $\sigma\bv=\{\sigma w \mid w\in \bv\}$.\label{item:systemdef2}
\item $(\integersnonnegative,\mathbb{W},\bv)$ is \emph{complete} if $\bv$ is closed in the topology of pointwise convergence. \label{item:systemdef3}
\end{enumerate}
\end{definition}
\smallskip
Note that if a dynamical system satisfies~\ref{item:systemdef1}-\ref{item:systemdef2} then~\ref{item:systemdef3} is equivalent to finite dimensionality of $\mathbb{W}$ (see~\cite[Section 7.1]{IM-JCW-SVH-BDM:06-book}). We denote the class of systems $(\integersnonnegative,\real^{m+p},\bv)$ satisfying~\ref{item:systemdef1}-\ref{item:systemdef3} by $\mathcal{L}^{m+p}$, where $m,p \in \integersnonnegative$. With slight abuse of notation and terminology, we denote a dynamical system in $\mathcal{L}^{m+p}$ only by its behaviour $\bv$.
\begin{definition}
The \emph{restricted behaviour in the interval $[1,T]$} is the set $\bv_T=\{w\in(\real^{m+p})^T \mid \exists\; v\in \bv\; \text{s.t. } w_t=v_t,\; 1\leq t\leq T\}$. A vector $w\in\bv_T$ is called a \emph{$T$-length trajectory of the dynamical system $\bv$}.
\end{definition}

\begin{definition}\label{def:controllable}
A system $\bv\in\mathcal{L}^{m+p}$ is \emph{controllable} if for every $T\in \integerspositive$, $w^1\in \bv_T$, $w^2\in\bv$ there exists $w\in \bv$ and $T'\in \integerspositive$ such that $w_t=w^1_t$ for $1\leq t \leq T$ and $w_t=w^2_{t-T-T'}$ for $t>T+T'$.
\end{definition}
\smallskip
In other words, a behavioural system is controllable if any two trajectories can be patched together in finite time. 

\subsection{Parametric system representation}\label{sec:parametric}

Without loss of generality, any trajectory $w\in \bv$ can be written as $w=\col(u,y)$, where $\col(u,y)\defeq (u^\top,y^\top)^\top$ (see~\cite[Theorem 2]{JCW:86}). In what follows, we will associate $u$ and $y$ with inputs and outputs.
There are several equivalent ways of representing a behavioural system $\bv \in \mathcal{L}^{m+p}$, including the classical input/output/state representation denoted by 
$\bv(A,B,C,D)=\{\col(u,y) \in (\real^{m+p})^{\integersnonnegative}\mid \exists\; x\in(\real^n)^{\integersnonnegative}\; \text{s.t.}\;  \sigma x=Ax+Bu,\; y=Cx+Du\}$. The input/output/state representation of smallest order (i.e., smallest state dimension) is called a \emph{minimal representation}, and we denote its order by $\bm{n}(\bv)$. Another important property of a system $\bv\in\mathcal{L}^{m+p}$ is the \emph{lag} defined by the smallest integer $\ell\in \integerspositive$ such that the observability matrix $\mathcal{O}_{\ell}(A,C)\defeq \col\left(C,CA, \dots,CA^{\ell-1}\right)$ has rank $\bm{n}(\bv)$. We denote the lag by $\bm{\ell}(\bv)$ (see~\cite[Section 7.2]{IM-JCW-SVH-BDM:06-book} for equivalent definitions). 
The lower triangular Toeplitz (impulse response) matrix consisting of $\Tf\in\integerspositive$ system Markov parameters is denoted by
\[
\mathcal{T}_\Tf(A,B,C,D) \defeq \begin{pmatrix}
D &0 &\cdots &0 \\
CB &D &\cdots &0 \\
\vdots &\ddots &\ddots &\vdots \\
CA^{\Tf-2}B &\cdots &CB &D
\end{pmatrix}.
\]
\begin{lemma}\longthmtitle{\hspace{1sp}\cite[Lemma 1]{IM-PR:08}}\label{lem:initialstate}
Let $\bv \in \mathcal{L}^{m+p}$ and $\bv(A,B,C,D)$ a minimal input/output/state representation. Let $\Tini,\Tf \in \integerspositive$ with $\Tini \geq \bm{\ell}(\bv)$ and $\col(\uini,u,\yini,y)\in \bv_{\Tini+\Tf}$. Then there exists a unique $\xini \in \real^{\bm{n}(\bv)}$ such that
\[
y=\mathcal{O}_\Tf(A,C)\xini +\mathcal{T}_\Tf(A,B,C,D)u.
\]
\end{lemma}
\smallskip
In other words, given a sufficiently long window of initial system data $\col(\uini,\yini)$, the initial state $\xini$ is unique and can be computed given knowledge of $A,B,C,D$ and $\col(u,y)$.

\subsection{Hankel and Page Matrices}
We introduce two important matrix structures: the Hankel matrix and the Page matrix.
\begin{definition}\label{def:pagehankel}
Let $L, T\in\integerspositive$. Let $u_{[1,T]}=\{u_j\}_{j=1}^T\subset \real^{m}$ be a sequence of vectors.
\begin{itemize}
\item We define the \emph{Hankel matrix\footnote{The classical definition of a Hankel matrix requires it to be square. We slightly abuse this classical terminology, allowing for general dimensions.} of depth $L$} as
\[
\mathscr{H}_L(u_{[1,T]})\defeq
\begin{pmatrix}
&u_1 &u_2 &\cdots &u_{T-L+1} \\
&u_2 &u_3 &\cdots &u_{T-L+2} \\
&\vdots &\vdots &\ddots &\vdots \\
&u_L &u_{L+1} &\cdots &u_T
\end{pmatrix}.
\]
\item We define the \emph{Page matrix of depth $L$} as
\[
\mathscr{P}_L(u_{[1,T]})\defeq
\begin{pmatrix}
&u_1 & u_{L+1} & \cdots & u_{(\lfloor \frac{T}{L} \rfloor -1)L+1} \\
&u_2 & u_{L+2} & \cdots & u_{(\lfloor \frac{T}{L} \rfloor -1)L+2} \\
& \vdots & \vdots & \ddots & \vdots \\
&u_L & u_{2L} & \cdots &u_{\lfloor \frac{T}{L} \rfloor L}
\end{pmatrix},
\]
where $\lfloor \cdot \rfloor$ is the floor function which rounds its argument down to the nearest integer.
\end{itemize}
\end{definition}
\smallskip
Note that if $T$ is a multiple of $L$, then $\lfloor \frac{T}{L} \rfloor L = T$ and the above expressions simplify accordingly. Hankel matrices have a long history in subspace identification~\cite{PVO-BDM:12-book} and more recently in data-driven control~\cite{IM-PR:08,JC-JL-FD:18,CDP-PT:19,JB-JK-MM-FA:20-journal}. Page matrices have also been used as an alternative to the classical Hankel matrix~\cite{AD-PVH-AH:82,AA-MJA-DS-DS:18}. 
\begin{definition}
Let $L,T,M\in\integerspositive$ and $u_{[1,T]}=\{u_j\}_{j=1}^T\subset \real^{m}$ be a sequence of vectors. We call $u_{[1,T]}$:
\begin{itemize}
\item \emph{Hankel exciting of order $L$} if the matrix $\mathscr{H}_L(u_{[1,T]})$ has full row rank.
\item \emph{$L$-Page exciting of order $M$} if the matrix 
\[
\begin{pmatrix}
\mathscr{P}_L(u_{[1,T-(M-1)L]}) \\
\mathscr{P}_L(u_{[L+1,T-(M-2)L]}) \\
\vdots \\
\mathscr{P}_L(u_{[L(M-1)+1,T]}) \\
\end{pmatrix}
\]
has full row rank.
\end{itemize}
\end{definition}
\smallskip
Roughly speaking, the terms Hankel and Page exciting refer to a collection of {inputs} sufficiently rich and long to yield an output sequence representative for the system's behaviour. Note that for a sequence of vectors to be Hankel exciting of order $L$, one requires that $T\geq L(m+1)-1$. For a sequence of vectors to be $L$-Page exciting of order $M$, one requires that $T\geq L((mL+1)M-1)$. We also observe that the definition of an $L$-Page exciting sequence of order $M$ depends on two indices, $L$ and $M$ (the number of block rows of the Page matrices and the number of vertically concatenated Page matrices, respectively), whereas the definition of a Hankel exciting sequence of order $L$ depends only on a single index $L$ (the number of block rows of the Hankel matrix).

The definition for Hankel exciting appeared in~\cite{JCW-PR-IM-BDM:05} under the name \emph{persistently exciting}. A more general notion of persistency of excitation appeared in~\cite{HJVW-CDP-MKC-PT:20} in the case when the data matrix was a mosaic Hankel matrix and was termed collective persistency of excitation. Our notion of Page exciting is a modification of these notions; the main difference is that there are no repeated entries in each of the data matrices. This has important implications when the entries of the matrix are corrupted by noise; see Section~\ref{sec:chooseparams}, the case-study \cite{LH-JC-JL-FD:20}, and references \cite{AD-PVH-AH:82,AA-MJA-DS-DS:18}. In short, the independence of the Page matrix entries leads to statistically and algorithmically favourable properties, e.g., singular-value thresholding can be used for de-noising. We will observe later that certain robustness and optimality guarantees become tight for Page matrices, and they lead to superior control performance. The price-to-pay is an increasing number of data samples. Almost all results within this paper hold for both Hankel and Page matrix structures, and we will explicitly comment on when the results differ depending on the {matrix structure chosen.}

\subsection{Non-parametric System Representation}
We now present a result known in behavioural systems theory as the \emph{Fundamental Lemma}~\cite{IM-PR:08}. The result first appeared in~\cite[Theorem 1]{JCW-PR-IM-BDM:05} in the case when the Hankel matrix was used, and was recently extended in~\cite[Theorem 2]{HJVW-CDP-MKC-PT:20} for more general mosaic Hankel matrices. We present the result using the Page matrix data structure introduced in Definition~\ref{def:pagehankel}.
\begin{theorem}\label{thm:fundamental}
Consider a controllable system $\bv\in\mathcal{L}^{m+p}$. Let $L,T\in\integerspositive$ with $L\geq \bm{n}(\bv)$. Let $(\hat{u}_{[1,T]},\hat{y}_{[1,T]})=\{(\hat{u}_j,\hat{y}_j)\}_{j=1}^T\subset \real^{(m+p)}$ be a $T$-length trajectory of $\bv$.  
Assume $\hat{u}_{[1,T]}$ to be $L$-Page exciting of order $\bm{n}(\bv)+1$. Then $(u_{[1,L]},y_{[1,L]})=\col(u_1,\dots,u_L,y_1,\dots,y_L)\in\bv_L$ if and only if there exists a vector $g\in\real^{\lfloor \frac{T}{L} \rfloor}$ such that
\begin{equation}\label{eq:datamatrixrep}
\begin{pmatrix}
\mathscr{P}_L(\hat{u}_{[1,T]})\\
\mathscr{P}_L(\hat{y}_{[1,T]})
\end{pmatrix}g
=\begin{pmatrix}
u_{[1,L]} \\
y_{[1,L]}
\end{pmatrix}
\end{equation}
\end{theorem}
\smallskip
The original result~\cite[Theorem 1]{JCW-PR-IM-BDM:05} required Hankel excitation of order $L+\bm{n}(\bv)$, and it coincides with Theorem~\ref{thm:fundamental} for $L=1$.
Theorem~\ref{thm:fundamental} replaces the need for a model or system identification process and allows for any trajectory of a controllable LTI system to be constructed using a finite number of data samples generated by a  sufficiently rich (i.e., $L$-page exciting) input sequence. Each column of the Page matrix in~\eqref{eq:datamatrixrep} is a trajectory of the system, and can be thought of as a motion primitive. By linearly combining elements of this trajectory library, we recover the whole space of trajectories.

In a sense, the Page matrix in~\eqref{eq:datamatrixrep} is a non-parametric predictive model based on raw data. To be precise in behavioural language, the Page matrix in~\eqref{eq:datamatrixrep} is an image representation of the restricted behaviour $\bv_L$. This is in contrast to kernel representations (with latent variables) such as state-space models which parameterize $\bv_L$ by means of its orthogonal complement. In what follows, we will express a linear system $\bv$ by its \emph{data matrix representation} in~\eqref{eq:datamatrixrep}. We refer to the Page matrix on the left-hand side of~\eqref{eq:datamatrixrep} as the {\em data matrix}.
\begin{remark}\label{rem:rank}
When $L\geq \bm{\ell}(\bv)$, it can be shown that the rank of the data matrix in the data matrix representation given by~\eqref{eq:datamatrixrep} is $mL+\bm{n}(\bv)$ where $m$ is the dimension of the inputs. Hence, the column span of the data matrix carves out a low-dimensional subspace of $\real^{L(m+p)}$ which coincides with the restricted behaviour $\bv_L$, the space of $L$-length trajectories. Since in general $\bm{n}(\bv)$ is unknown, it suffices to replace $\bm{n}(\bv)$ in Theorem~\ref{thm:fundamental} with an upper bound on the system order.
\oprocend
\end{remark}
For the rest of the paper, we think of an LTI system $\bv$ only in terms of its data matrix representation. We see next that Theorem~\ref{thm:fundamental} allows to implicitly estimate the state, predict the future behaviour, and design optimal control inputs~\cite{IM-PR:08}.
\section{Deterministic DeePC}\label{sec:detdeepc}

Consider the controllable LTI system $\bv\in\mathcal{L}^{m+p}$ whose model is unknown. 
We address the problem of finite-horizon optimal control, where the goal is to design a finite sequence of control inputs that result in desirable outputs. In particular, let $t\in\integersnonnegative$ be the current time, let $\Tf\in\integerspositive$ be the prediction horizon, and let $f_1\colon \real^{m\Tf}\to \realnonnegative$ (respectively, $f_2\colon \real^{p\Tf}\to \realnonnegative$) be a cost function on the future inputs (respectively, outputs). We wish to design an input sequence $\col(u_t,\dots,u_{t+\Tf-1})\in\real^{m\Tf}$ such that the corresponding output sequence $\col(y_t,\dots,y_{t+\Tf-1})\in\real^{p\Tf}$ minimizes the cost $f_1+f_2$. Furthermore, we wish for the inputs and outputs to lie in the constraint sets $\mathcal{U}\subseteq\real^{m\Tf}$ and $\mathcal{Y}\subseteq \real^{p\Tf}$, respectively.

The conventional formulation of this finite-time optimal control problem uses a parametric input/output/state representation $\bv(A,B,C,D)$ for prediction and estimation:
\begin{equation}\label{eq:mpc}
\begin{aligned}
\underset{u,y,x}{\min}\quad
& f_1(u)+f_2(y) \\
\text{s.t.\quad}
& x_{k+1}=Ax_k+Bu_k \; \forall k\in\{0,\dots,\Tf-1\} \\
&y_k=Cx_k+Du_k \; \forall k\in\{0,\dots,\Tf-1\} \\
&x_0 = \hat{x}_t \\
&u \in\mathcal{U},\; y \in\mathcal{Y}\,.
\end{aligned}
\end{equation}
Here, $\hat{x}_t$ is the current state at time $t$, which in the deterministic case can be obtained, e.g., by propagating the model equations backward in time by $\Tini\geq \bm{\ell}(\bv)$ steps; see Lemma~\ref{lem:initialstate}. When optimization problem~\eqref{eq:mpc} is implemented in a receding horizon fashion it is known as model-predictive control (MPC) and is a widely celebrated control technique. Typically, in receding horizon MPC problems the constraints take the form of separable stage constraints with an additional terminal constraint (similarly for the cost functions)~\cite{FB-AB-MM:17-book}.

Using Theorem~\ref{thm:fundamental}, we can see that raw data collected from system $\bv$ can be used to perform implicit state estimation at the same time as predict forward trajectories of $\bv$. Indeed, let $\Tini,\Tf\in\integerspositive$ and let $(\hat{u}_{[1,T]},\hat{y}_{[1,T]})=\{(\hat{u}_j,\hat{y}_j)\}_{j=1}^T\subset \real^{(m+p)}$ be a $T$-length trajectory of $\bv$ collected offline. Assume that $\hat{u}_{[1,T]}$ is $(\Tini+\Tf)$-Page exciting of order $\bm{n}(\bv)+1$. We partition the data matrices into two parts; one will be used to perform the implicit state estimation and the other to perform forward predictions. More formally, we use the language of subspace system identification~\cite{PVO-BDM:12-book} and define
\begin{equation}\label{eq:data}
\begin{pmatrix}
\Uphat\\
\Ufhat
\end{pmatrix}
\defeq \mathscr{P}_{\Tini+\Tf}(\hat{u}_{[1,T]}),\quad
\begin{pmatrix}
\Yphat\\
\Yfhat
\end{pmatrix}
\defeq \mathscr{P}_{\Tini+\Tf}(\hat{y}_{[1,T]}),
\end{equation}
where $\Uphat\in\real^{m\Tini\times {\lfloor \frac{T}{\Tini+\Tf} \rfloor}}$ consists of the first $\Tini$ block rows of $\mathscr{P}_{\Tini+\Tf}(\hat{u}_{[1,T]})$ and $\Ufhat\in\real^{m\Tf\times {\lfloor \frac{T}{\Tini+\Tf} \rfloor}}$ the last $\Tf$ block rows (similarly for $\Yphat$ and $\Yfhat$). Let $t\in\integersnonnegative$ be the current time and denote by $\col(\uinihat,\yinihat)$ the most recent $\Tini$-length trajectory measured from $\bv$. Then given a sequence of inputs $u=\col(u_t,\dots,u_{\Tf-1})$ of length $\Tf\in\integerspositive$, one may solve for a vector $g\in\real^{\lfloor \frac{T}{\Tini+\Tf} \rfloor}$ satisfying
\begin{equation}\label{eq:predict}
\begin{pmatrix}
\Uphat\\
\Yphat\\
\Ufhat
\end{pmatrix}
g = \begin{pmatrix}
\hat{u}_{\textup{ini}} \\
\hat{y}_{\textup{ini}} \\
u
\end{pmatrix}.
\end{equation}
By Theorem~\ref{thm:fundamental}, the predicted $\Tf$-length output of $\bv$ is then given by $y=\Yfhat g$. Note that in general the $g$ satisfying~\eqref{eq:predict} is non-unique. However, by Lemma~\ref{lem:initialstate}, if $\Tini\geq \bm{\ell}(\bv)$ then the predicted output $y$ is unique. The top two block equations of~\eqref{eq:predict} can be thought of as implicitly fixing the initial state from which the future trajectory will depart. 

We now recall the deterministic DeePC method presented in~\cite{JC-JL-FD:18} (using Hankel matrices) for controlling the system $\bv$.

\noindent\underline{\textbf{Data Collection} (offline):}
Apply to $\bv$ a $T$-length input sequence $\hat{u}_{[1,T]}$ that is $(\Tini+\Tf)$-Page exciting of order $\bm{n}(\bv)+1$. Measure the corresponding output trajectory $\hat{y}_{[1,T]}$. Form data matrices $\Uphat$, $\Ufhat$, $\Yphat$, $\Yfhat$ according to~\eqref{eq:data}.

\noindent\underline{\textbf{DeePC} (online):}
Let $f_1\colon \real^{m\Tf}\to\realnonnegative$, $f_2\colon\real^{p\Tf}\to\realnonnegative$ be cost functions describing the costs on future inputs and outputs of $\bv$, respectively. Collect the most recent $\Tini$-length trajectory $\col(\uinihat,\yinihat)$ measured from $\bv$. Solve the finite-horizon optimal control problem:
\begin{equation}\label{eq:deepc}
\begin{aligned}
\underset{g}{\min}\quad
& f_1(\Ufhat g)+f_2(\Yfhat g) \\
\text{s.t.\quad}
& \begin{pmatrix}
\Uphat \\ \Yphat
\end{pmatrix}g
=\begin{pmatrix}
\hat{u}_{\textup{ini}} \\ 
\hat{y}_{\textup{ini}}
\end{pmatrix}\\
&\Ufhat g \in\mathcal{U}\\
&\Yfhat g \in\mathcal{Y}.
\end{aligned}
\end{equation}
Similar to MPC, the online DeePC optimization can be implemented in a receding horizon fashion.
Note that the optimization problem~\eqref{eq:deepc} only requires input/output measurements from the system and does not require the identification of a model. The equivalence of DeePC to classical model predictive control when a system model for $\bv$ is given, is established in~\cite{JC-JL-FD:18} in the case when the data matrix is a Hankel matrix.
\begin{proposition}\label{thm:setequiv}
Consider an LTI system $\bv$ whose input/output/state representation is given by $\bv(A,B,C,D)$. Consider the MPC optimization problem \eqref{eq:mpc}
 and the optimization problem~\eqref{eq:deepc} with $\Tini\geq \bm{\ell}(\bv)$. Then $g$ satisfying the constraints of~\eqref{eq:deepc} implies that $(u,y)=(\Ufhat g,\Yfhat g)$ satisfies the constraints of~\eqref{eq:mpc}. Conversely, $(u,y)$ satisfying the constraints of~\eqref{eq:mpc} implies the existence of $g$ satisfying the constraints of~\eqref{eq:deepc} where $(u,y) = (\Ufhat g,\Yfhat g)$.
\end{proposition}
\smallskip
As a direct corollary of Proposition~\ref{thm:setequiv}, DeePC problem~\eqref{eq:deepc} and MPC problem~\eqref{eq:mpc} result in equivalent and unique (assuming convexity) open-loop (resp., closed-loop) behaviour when applied in an open-loop (resp., receding horizon) fashion to a deterministic LTI system $\bv$. As a major difference, most MPC approaches assume direct state measurements, which allows certain methods to be applied (e.g., inclusion of terminal state constraints). We remark that $(\uini,\yini,u,y)$ are coordinates of a (generally non-minimal) realization allowing for similar methodologies. See~\cite{JB-JK-MM-FA:20-journal} in which stability and recursive feasibility for a variant of~\eqref{eq:deepc} have been shown.

\section{Distributionally Robust DeePC}\label{sec:stodeepc}

We now consider the case where system $\bv$ is subject to stochastic disturbances, which result in noisy output measurements in the data matrices. We begin with modifying the objective function and constraints of the deterministic DeePC problem~\eqref{eq:deepc} to robustify against the stochastic disturbances. This leads to a semi-infinite optimization problem which we reformulate under reasonable assumptions into a finite, convex, and tractable program termed the distributionally robust DeePC problem. Finally, we show that the distributionally robust DeePC problem enjoys robust performance guarantees.

\subsection{Problem Setup}

We view the output data matrices $\Yphat$ and $\Yfhat$ given by~\eqref{eq:data} as particular realizations of random variables which we denote by $\Yp$ and $\Yf$, respectively. Recall that we use the $\widehat{\cdot}$ notation to denote recorded data samples and to indicate that objects depend on recorded data samples. Define the random variable 
$\xi=(\xi_1^{\top},\dots,\xi_{p(\Tini+\Tf)}^\top)^\top$ where $\xi_i^\top\in\real^{\lfloor \frac{T}{\Tini+\Tf}\rfloor}$ denotes the $i$-th row of the matrix $(\Yp^\top\;\Yf^\top)^\top$. We denote the support set of the random variable $\xi$ by $\Xi$. We denote the probability distribution of the random variable $\xi$ by $\Prob$. Note that this distribution is induced by the system $\bv$ itself, and any stochastic disturbance acting on the system. 

Due to the stochastic nature of the system, the future trajectory predictions obtained from the deterministic DeePC optimization problem~\eqref{eq:deepc} are uncertain. Furthermore, the realized data matrix~\eqref{eq:datamatrixrep} may not describe a low-dimensional subspace as described in Remark~\ref{rem:rank}. In fact, the data matrix in~\eqref{eq:datamatrixrep} will likely have full rank and can thus predict arbitrary future trajectories. Finally, the distribution $\Prob$ itself is  unknown since we do not have models of the system and disturbances. Hence, we make some changes to the objective function and the constraints in~\eqref{eq:deepc} to robustify against the noisy data.

\emph{Objective Function:} Since the constraint used to obtain the initial condition of the future trajectory given by $\Yphat g=\yinihat$ is affected by the stochastic disturbance and may not be feasible, we lift it into the objective function using an estimation function $f_3\colon \real^{p\Tini}\to\realnonnegative$ mapping $\Yp g -\yinihat$ to a penalty. Including $f_3$ in the objective function can be thought of as an implicit estimation of the initial condition from which the predicted future trajectory must evolve. Choosing $f_3(\cdot)=\|\cdot\|_2$ would result in a least-squares type initial condition estimate reminiscent of moving horizon estimation~\cite{FB-AB-MM:17-book}.

Since the system is stochastic, we focus on minimizing the expectation of the objective function in~\eqref{eq:deepc}, where the expectation is taken with respect to the true distribution $\Prob$, i.e., we wish to minimize $\Exp_{\Prob}[f_1(\Ufhat g)+f_2(\Yf g)+f_3(\Yp g-\yinihat)]$. 

Recall that the distribution $\Prob$ pertains to the offline samples $(\hat{u}_{[1,T]},\hat{y}_{[1,T]})$ and not to the real-time measurement $\yinihat$ addressed through the above least-square estimation. As discussed before, the distribution $\Prob$ itself is unknown, and must be estimated from data. In order to be robust to errors in this estimate, we use distributionally robust optimization techniques; namely, we construct a so-called ambiguity set $\ambiguity$ depending  on the collected data (specified precisely later) such that the true distribution lies in the ambiguity set with high confidence.  We then optimize the expectation of the cost function, where the expectation is taken with respect to the worst-case distribution in $\ambiguity$. This leads us to the distributionally robust objective
\[
\underset{\Q\in\ambiguity}{\sup} \Exp_{\Q}\left[f_1(\Ufhat g) + f_2(\Yf g) + f_3(\Yp g -\yinihat) \right].
\]
We note that this formulation leads to robustness against a ``set of systems'' compatible with the (possibly noisy) data samples, and that this ``set of systems'' is broader than mere LTI systems with additive process and measurement noise.

\begin{figure}[h]
\begin{center}
\includegraphics[width=0.45\textwidth]{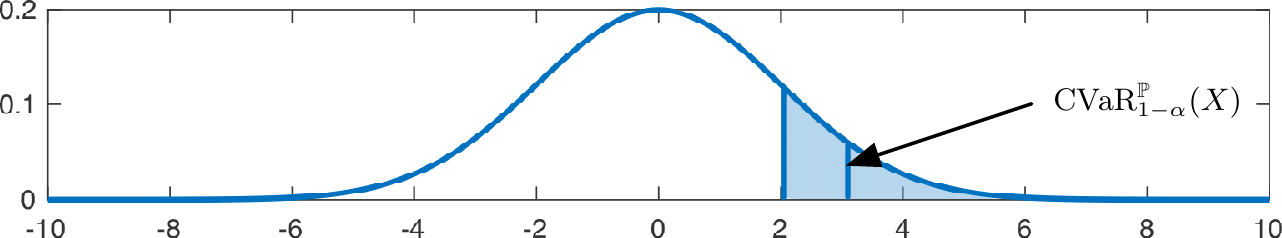}
\end{center}
\caption{Depiction of the conditional value at risk at level $\alpha=0.15$ for a Gaussian random variable $X$. The shaded region accounts for $\alpha \times 100 \%$ of the mass of the Gaussian distribution.}
\label{fig:cvar}
\end{figure}

\emph{Constraints:} Requiring the future outputs to lie in a particular subset $\mathcal{Y}$ almost surely as in~\eqref{eq:deepc} may not be possible. Thus, we relax the hard output constraint to a so-called conditional value at risk (CVaR) constraint. Let the output constraint be written as $\mathcal{Y}=\{y\;\vert\; h(y)\leq 0 \}$, where $h\colon \real^{p\Tf}\to \real$ describes desired constraints on future trajectories $y$. We define
\[
\textup{CVaR}_{1-\alpha}^{\Prob}(h(y))\defeq \inf_{\tau\in\real}\left\{\tau+\frac{1}{\alpha}\Exp_{\Prob}[(h(y)-\tau)_{+}]\right\},
\]
where $\alpha\in(0,1)$ is a user-chosen confidence parameter. It is well known that the constraint $\textup{CVaR}_{1-\alpha}^{\Prob}(h(y))\leq 0$ is a convex relaxation of the classical chance constraint (or VaR) given by $\Prob(h(y)\leq 0)\geq 1-\alpha$ (see, e.g., \cite{AN-AS:06}). In other words, the CVaR constraint ensures that the constraint $y\in\mathcal{Y}$ will hold with high probability. In the case when the stochastic disturbance is a continuous random variable, the CVaR constraint can be interpreted as penalizing the expected violation in the $\alpha$ percent of cases where violations do occur (see Figure~\ref{fig:cvar}). The latter point is particularly important in control applications where large violations of constraints could lead to catastrophic system behaviour. Hence, the CVaR constraint is not merely a relaxation of the classical VaR constraint, but often the more reasonable constraint formulation for control problems.

\noindent Similar to minimizing the worst case expected objective function, we ask that the CVaR constraint be satisfied for the worst-case distribution in $\ambiguity$. This leads us to the robust constraint
\[
\underset{\Q\in\ambiguity}{\sup}\textup{CVaR}_{1-\alpha}^{\Q}(h(y))\leq 0.
\]

One particular ambiguity set $\ambiguity$ that offers strong performance guarantees as well as mathematical tractability while allowing the user to adjust its conservatism is the \emph{Wasserstein ambiguity set} centred around the so-called \emph{empirical distribution}, i.e., the sample distribution~\cite{PME-DK:18}. In particular, let $\mathcal{M}(\Xi)$ be the set of all distributions $\Q$ supported on $\Xi$ such that $\Exp_\Q[\|\xi\|_{r}]<\infty$, where $\|\cdot\|_{r}$ is the $r$-norm for some $r\in[1,\infty]$.
\begin{definition}
Let $r\in[1,\infty]$. The \emph{Wasserstein metric} $\textup{d}_\textup{W}\colon\mathcal{M}(\Xi)\times\mathcal{M}(\Xi)\to\realnonnegative$ is defined as
\[
\textup{d}_\textup{W}(\Q_1,\Q_2)\defeq\inf_{\Pi}\left\{\int_{\Xi^2}\|\xi_1-\xi_2\|_{r}\Pi(\textup{d}\xi_1,\textup{d}\xi_2)\right\},
\]
where $\Pi$ is a joint distribution of $\xi_1$ and $\xi_2$ with marginal distributions $\Q_1\in \mathcal{M}(\Xi)$ and $\Q_2\in\mathcal{M}(\Xi)$, respectively.
\end{definition}
\smallskip
The Wasserstein metric can be viewed as a distance between probability distributions, where the distance is calculated via an optimal mass transport plan $\Pi$.
For $\eps\geq0$, we denote the \emph{Wasserstein ball of radius $\eps$ centred around distribution $\Q$} by 
\[
B_{\eps}(\Q)\defeq\{\Q'\in\mathcal{M}(\Xi)\;\vert\; \textup{d}_{\textup{W}}(\Q,\Q')\leq \eps\}\,.
\]
We chose the Wasserstein ball as the ambiguity set as opposed to other popular ambiguity sets such as $f$-divergence or moment-based ambiguity sets since the latter require prior knowledge and assumptions on the true data generating distribution (absolute continuity with respect to a nominal distribution or known moments).

We now present the distributionally robust DeePC method for controlling stochastic systems.

\noindent\underline{\textbf{Data Collection} (offline):} 
We collect $N\in\integerspositive$ trajectories via repeated identical experiments. We start by fixing an input sequence $\hat{u}_{[1,T]}$ that is $(\Tini+\Tf)$-Page exciting of order $\bm{n}(\bv)+1$. For each experiment $i\in\{1,\dots,N\}$, we assume that the system is initialized to the same (unknown) state, apply inputs $\hat{u}_{[1,T]}$, and measure the corresponding output trajectory $\hat{y}_{[1,T]}^{(i)}$. The superscript $(i)$ denotes the trajectory obtained from the $i$-th experiment. For each data batch $i\in\{1,\dots,N\}$, build data matrices $\Uphat,\Ufhat$, $\Yphat^{(i)}$, and $\Yfhat^{(i)}$ using~\eqref{eq:data}. The output data matrices implicitly define $N$ data samples which we denote by $\hat{\xi}^{(i)}$ for $i\in\{1,\dots,N\}$. Lastly, we use the data collected to build the \emph{empirical distribution} $\Phat_N=\frac{1}{N}\sum_{i=1}^N \delta_{\hat{\xi}^{(i)}}$. 
Note that collecting only one $T$-length trajectory (i.e., $N=1$) is possible. As we will see in Section~\ref{sec:chooseparams} though, choosing $N$ larger can lead to less uncertainty and better performance.

\noindent\underline{\textbf{Distributionally Robust DeePC} (online):} We consider the Wasserstein ball of radius $\eps$ (which will be chosen as a hyperparameter quantifying a desired robustness level) around the empirical distribution $\Phat_N$ as the ambiguity set. Hence, the distributionally robust DeePC problem is given by
\begin{equation}\label{eq:robustdeepc}
\begin{aligned}
\underset{g}{\min}\;\underset{\Q\in\B}{\sup} \quad & \Exp_{\Q}\left[f_1(\Ufhat g) + f_2(\Yf g) + f_3(\Yp g -\yinihat) \right] \\
\text{s.t.\quad}
&\Uphat g = \hat{u}_{\textup{ini}}\\
&\Ufhat g \in\mathcal{U}\\
&\underset{\Q\in\B}{\sup}\textup{CVaR}_{1-\alpha}^{\Q}(h(\Yf g))\leq 0
\end{aligned}
\end{equation}

By assuming that the system is initialized to the same (unknown) state at the beginning of each experiment and applying an identical input sequence for each experiment, we ensure that the collection of data matrices $\Yphat^{(i)}$ and $\Yfhat^{(i)}$ describe $N$ independent and identically distributed (i.i.d.) samples of the random variable $\xi$. Before their realization, the collection $\{\hat{\xi}^{(i)}\}_{i=1}^N$ can be viewed as a random object governed by the $N$-fold product distribution $\Prob^N$. The empirical distribution $\Phat_N$ is our sample average estimate for the true distribution $\Prob$ and is set as the centre of the ambiguity set in~\eqref{eq:robustdeepc}.

Compared to MPC or deterministic DeePC discussed in Section~\ref{sec:detdeepc}, the distributional robust program \eqref{eq:robustdeepc} features the additional  robustness parameter $\eps$, whose selection will be discussed later, and $\alpha$ whose value is usually chosen to be 0.1, 0.05, or 0.01. One notable feature of~\eqref{eq:robustdeepc} is that it is robust to all probability distributions within the ambiguity set that could describe the data matrices $\Yp$ and $\Yf$. Hence, the above formulation is robust to a ``set of systems'' which captures more than LTI systems with additive noise. In fact, for a sufficiently large robustness parameter $\eps$, the ambiguity set considers all time series originating from linear or nonlinear systems as long as they have ``finite expectation'', i.e., the integrals defining the Wasserstein metric exist. This follows from the fact that the Wasserstein ball contains all Dirac distributions generated by time series $\xi$ satisfying $\textup{d}_{\textup{W}}(\delta_{\xi},\Phat_N)\leq\eps$. Hence, for $\eps$ large enough, the time series $\xi$ can be arbitrary.

Problem~\eqref{eq:robustdeepc} is a semi-infinite optimization problem due to the supremums over the space of distributions in the cost and in the constraints.
Moreover, the Wasserstein metric (resp. CVaR) are themselves defined via infinite (resp. finite) programs. Thus, \eqref{eq:robustdeepc} does not immediately present itself as a tractable program. Below, we show that under reasonable assumptions on the objective and constraint function, the semi-infinite optimization problem above admits a tractable reformulation. 

\subsection{Main results}
We provide a tractable reformulation of the distributionally robust DeePC problem~\eqref{eq:robustdeepc} which depends on two intermediate lemmas (Lemma~\ref{lem:costreformulation} and~\ref{lem:constraintreformulation_lipschitz} in the Appendix) providing separate reformulations for the objective function and constraints.
In the statement of the result below, we make use of the dual norm $\|\cdot\|_{q} = \|\cdot\|_{r,\ast}$, where, $q$ and $r$ satisfy $\frac{1}{r}+\frac{1}{q}=1$.
\begin{theorem}\longthmtitle{Tractable Reformulation}\label{thm:reformulation}
Assume that $f_1$ is convex and $f_2$, $f_3$, and $h$ are convex and Lipschitz continuous. Specifically, let $L_{\textup{obj}}>0$ (resp., $L_{\textup{con}}>0$) be the Lipschitz constant with respect to the $r$-norm of the  mapping $(x,y)\mapsto f_2(x)+f_3(y)$ (resp., $h$). Then the optimal value of problem~\eqref{eq:robustdeepc} is upper bounded by the optimal value of
\begin{equation}\label{eq:tractabledeepc}
\begin{aligned}
\underset{g,\tau,s_i}{\min} \quad &f_1(\Ufhat g)+\frac{1}{N}\sum_{i=1}^N \left(f_2(\Yfhat^{(i)}g)+f_3(\Yphat^{(i)} g-\yinihat)\right)\\
& \quad+L_{\textup{obj}}\eps \|g\|_{r,\ast}\\
\text{s.t.\quad} &\Uphat g=\uinihat \\
&\Ufhat g\in\mathcal{U} \\
&-\tau \alpha + L_{\textup{con}}\eps \|g\|_{r,\ast} +\frac{1}{N}\sum_{i=1}^N s_i\leq 0\\
&\tau + h(\Yfhat^{(i)} g) \leq s_i \quad \forall i\leq N\\
&s_i\geq 0 \quad \forall i\leq N.
\end{aligned}
\end{equation}
Moreover, the solution $\hat{g}^{\star}$ to the above will satisfy $\textup{CVaR}_{1-\alpha}^{\Q}(h(\Yf \hat{g}^{\star}))\leq 0$ for all $\Q\in\B$. The objective of~\eqref{eq:tractabledeepc} coincides with the objective of~\eqref{eq:robustdeepc} when $\Xi=\real^{p(\Tini+\Tf){\lfloor \frac{T}{\Tini+\Tf}\rfloor}}$.
\end{theorem}
\smallskip
Observe that when the data matrices in~\eqref{eq:data} are Hankel matrices, then $\Xi$ is always a strict subspace of $\real^{p(\Tini+\Tf){\lfloor \frac{T}{\Tini+\Tf}\rfloor}}$ encoding the Hankel structure. When the matrices are Page matrices, the support set $\Xi$ may coincide with $\real^{p(\Tini+\Tf){\lfloor \frac{T}{\Tini+\Tf}\rfloor}}$ (e.g., if the measurements are affected by Gaussian noise).
A similar result holds for piecewise affine constraint functions by combining Lemma~\ref{lem:costreformulation} and~\ref{lem:constraintreformulation_affine}. In this case, the reformulated constraint set can be made tight at the cost of making the objective an upperbound (see Remark~\ref{rem:tight}). Since, we treat the objective and constraints separately (via Lemma~\ref{lem:costreformulation} and~\ref{lem:constraintreformulation_lipschitz}), it is possible to choose different values of $\eps$ in the constraint and objective of~\eqref{eq:robustdeepc}, though we refrain from doing so for clarity of exposition. Furthermore, by treating the objective and the constraints separately, their worst case distributions may differ leading to a more conservative solution.

The above result shows that the distributionally robust objectives can be reformulated as the sample average of the objective function plus an additional regularization term. To the best of our knowledge this observation was first made in~\cite{SSA-DK-PME:19} in the context of machine learning problems. Hence, being robust in the trajectory space (data space) in the sense of the $r$-norm requires regularizing the sample average objective with the \emph{dual} $q$-norm, where the weight of regularization depends on the Lipschitz constant of the objective function and the radius of the Wasserstein ball. For example, the 1-norm regularization adopted in \cite{JC-JL-FD:18} corresponds to $\infty$-norm robustness in the trajectory space.

We now show that the robust DeePC problem~\eqref{eq:tractabledeepc} exhibits strong probabilistic guarantees. We need the following auxiliary measure concentration result under a light-tailedness assumption.
\begin{theorem}\label{thm:measureconcentration}\longthmtitle{Measure Concentration~\cite[Theorem 2]{NF-AG:15}}
Assume that $\Exp_{\Prob}[\textup{exp}(\|\xi\|_r^a)]<\infty$ for some $a>1$. Then
\[
\begin{aligned}
&\Prob^N\left\{d_{\textup{W}}(\Prob,\Phat_N)\geq \eps\right\}\\
&\leq \begin{cases}
c_1 \textup{exp}(-c_2 N \eps^{\lfloor \frac{T}{\Tini+\Tf}\rfloor}) & \textup{if }\eps\leq 1\\
c_1 \textup{exp}(-c_2 N \eps^a) & \textup{if }\eps> 1
\end{cases}
\end{aligned}
\]
for all $N\geq 1$, and $\eps>0$, where $\Prob^N$ is the $N$-fold product distribution, and $c_1, c_2$ are positive constants depending on $a$, ${\lfloor \frac{T}{\Tini+\Tf}\rfloor}$ and the value of $\Exp_{\Prob}[\textup{exp}(\|\xi\|_r^a)]$.
\end{theorem}
\smallskip

We can use Theorem~\ref{thm:measureconcentration} in order to compute the minimum Wasserstein radius $\eps$ such that the true data-generating distribution $\Prob$ lives inside the Wasserstein ball $\B$ with confidence at least $1-\beta$. Indeed, by inverting the above inequalities the minimum radius is given by
\begin{equation}\label{eq:mineps}
\eps(\beta,N) = \begin{cases}
\left(\frac{\log(c_1\beta^{-1})}{c_2 N}\right)^{\frac{1}{\lfloor \frac{T}{\Tini+\Tf}\rfloor}} & \textup{if } N\geq \frac{\log(c_1\beta^{-1})}{c_2}\\
\left(\frac{\log(c_1\beta^{-1})}{c_2 N}\right)^{\frac{1}{a}} & \textup{if } N< \frac{\log(c_1\beta^{-1})}{c_2}
\end{cases}
\end{equation}
\begin{theorem}\label{thm:performanceguarantee}\longthmtitle{Robust Performance Guarantee}
Assume that  $\Exp_{\Prob}[e^{\|\xi\|_{r}^a}]<\infty$ for some $a>1$. Let $\beta\in(0,1)$ be the desired level of confidence. Let $J(g)\defeq \Exp_{\Prob}[ f_1(\Ufhat g) +f_2(\Yf g)+f_3(\Yp g -\yinihat)]$. Let $\widehat{J}(g)$ denote the value of the objective function in~\eqref{eq:tractabledeepc} evaluated at $g$ and $\hat{g}^{\star}$ denote an optimizer of problem~\eqref{eq:tractabledeepc} with $\eps(\beta,N)$ chosen as in~\eqref{eq:mineps}.  Then
\[
\begin{aligned}
&\Prob^N\left\{J(\hat{g}^{\star})\leq \widehat{J}(\hat{g}^{\star})\right\}\geq 1-\beta,\;\textup{and} \\
&\Prob^N\left\{\textup{CVaR}_{1-\alpha}^{\Prob}(h(\Yf \hat{g}^{\star}))\leq 0\right\}\geq 1-\beta.
\end{aligned}
\]
\end{theorem}
\smallskip

Theorem~\ref{thm:performanceguarantee} shows that with probability $1-\beta$ (with respect to the $N$-fold product distribution $\Prob^N$ of $\Prob$), the optimal value of the robust DeePC problem~\eqref{eq:tractabledeepc} is an upper bound for the expected cost with respect to the true distribution $\Prob$. Furthermore, the CVaR constraint  with respect to the true distribution $\Prob$ holds with probability $1-\beta$.
\begin{remark}\label{rem:radius}
The Wasserstein radius given by~\eqref{eq:mineps} decays with a rate of $N^{\frac{1}{\lfloor T/(\Tini+\Tf)\rfloor}}$. To decrease the radius by 50\%, the number of samples $N$ must increase by $2^{\lfloor \frac{T}{\Tini+\Tf}\rfloor}$. This rate (albeit unfortunate) is tight~\cite{DK-PME-VN-SS:19}. In practice, the Wasserstein ball radius given by~\eqref{eq:mineps} is often larger than necessary, i.e., $\Prob\not\in \B$ with probability much less than $\beta$. Furthermore, even when $\Prob\not\in \B$, the robust quantity $\widehat{J}(\hat{g}^{\star})$ may still serve as an upper bound for the out-of-sample performance $J(\hat{g}^{\star})$~\cite{PME-DK:18}. In addition to this conservativeness, we note that the Wasserstein radius provided by~\eqref{eq:mineps} is hardly quantifiable as it depends on constants which are unknown to us in our data-driven setting. For practical purposes, one should choose the radius of $\B$ in a data-driven fashion (e.g., as we have done via tuning in experiments in Section~\ref{sec:chooseparams}). Alternatively, one could use standard cross-validation techniques (some of which are stated in [5]) or more advanced techniques (albeit only providing asymptotic guarantees) provided in~\cite{JB-YK-KM:19}. \oprocend
\end{remark}

The robust DeePC Algorithm \ref{alg:robustdeepc} implements~\eqref{eq:tractabledeepc} in a receding horizon fashion. 
We note that all results in this section hold when the data matrices are Hankel matrices in~\eqref{eq:data}. However, the upper bound and set inclusions presented in Lemmas~\ref{lem:costreformulation},~\ref{lem:constraintreformulation_lipschitz} and Theorem~\ref{thm:reformulation} are not tight because the Hankel structure restricts $\Xi$ to a strict subspace of $\real^{p(\Tini+\Tf)\lfloor \frac{T}{\Tini+\Tf}\rfloor}$. If the constraints are piecewise affine (see Lemma~\ref{lem:constraintreformulation_affine}), then the constraint reformulation can be tight for Hankel matrices.
In the following case study we will compare performance of Hankel and Page matrix structures.
\setlength{\intextsep}{10pt}
\begin{algorithm}[h]
	\caption{Robust DeePC}
	\label{alg:robustdeepc}
	\textbf{Input:} 
	data trajectories $\col(\hat{u},\hat{y}^{(i)})\in\real^{(m+p)T}$ for $i\in\{1,\dots,N\}$ with $\hat{u}$ $(\Tini+\Tf)$-Page exciting of order $\bm{n}(\bv)+1$,\\
	 most recent input/output measurements $\col(\hat{u}_{\textup{ini}},\hat{y}_{\textup{ini}})\in\real^{(m+p)\Tini}$
	\begin{algorithmic}[1]
		\STATE \label{step:deepcbegin} Solve~\eqref{eq:tractabledeepc} for $\hat{g}^{\star}$.
		\STATE Compute optimal input sequence $u^{\star} = \Ufhat \hat{g}^{\star}$.
		\STATE Apply optimal input sequence $(u_t,\dots,u_{t+\nu-1})=(u_1^{\star},\dots,u_{\nu}^{\star})$ for some $\nu\leq \Tf$.
		\STATE Set $t$ to $t+\nu$ and update $\hat{u}_{\textup{ini}}$ and $\hat{y}_{\textup{ini}}$ to the $\Tini$ most recent input/output measurements.
		\STATE Return to~\ref{step:deepcbegin}.
	\end{algorithmic}
\end{algorithm}
%

\section{Numerical Results}\label{sec:numericalresults}

\subsection{Nonlinear and Stochastic Aerial Robotics Case Study}\label{sec:casestudy}
We illustrate the performance of distributionally robust DeePC with a high-fidelity nonlinear quadcopter simulation~\cite{EE-JC-PB-JL-FD:19}. The output measurements are the 3 spatial coordinates $(p_{x},p_{y},p_{z})$ of the quadcopter. We denote the output at time $t$ by $y_t=(p_{x,t},p_{y,t},p_{z,t})\in \real^3$ and the $i$-th component of $y_t$ by $y_{t,i}$. The inputs are the total thrust produced by all 4 rotors $f_{\textup{tot}}$, and the angular body rates around the $x$ and $y$ body axes of the quadcopter, $\omega_{x},\omega_{y}$ respectively. We denote the input at time $t$ by $u_t=(f_{\textup{tot},t},\omega_{x,t},\omega_{y,t})\in\real^3$ and the $i$-th component of $u_t$ by $u_{t,i}$. 
The output measurements are affected by additive zero-mean Gaussian noise during the data collection phase (offline) and the control phase (online) in which the distributionally robust DeePC Algorithm~\ref{alg:robustdeepc} is implemented. Statistics of the noise were chosen to closely match the experimental setup~\cite{EE-JC-PB-JL-FD:19}.

\begin{figure*}[b]
    \centering 
\begin{subfigure}{0.33\linewidth}
  \includegraphics[width=\linewidth]{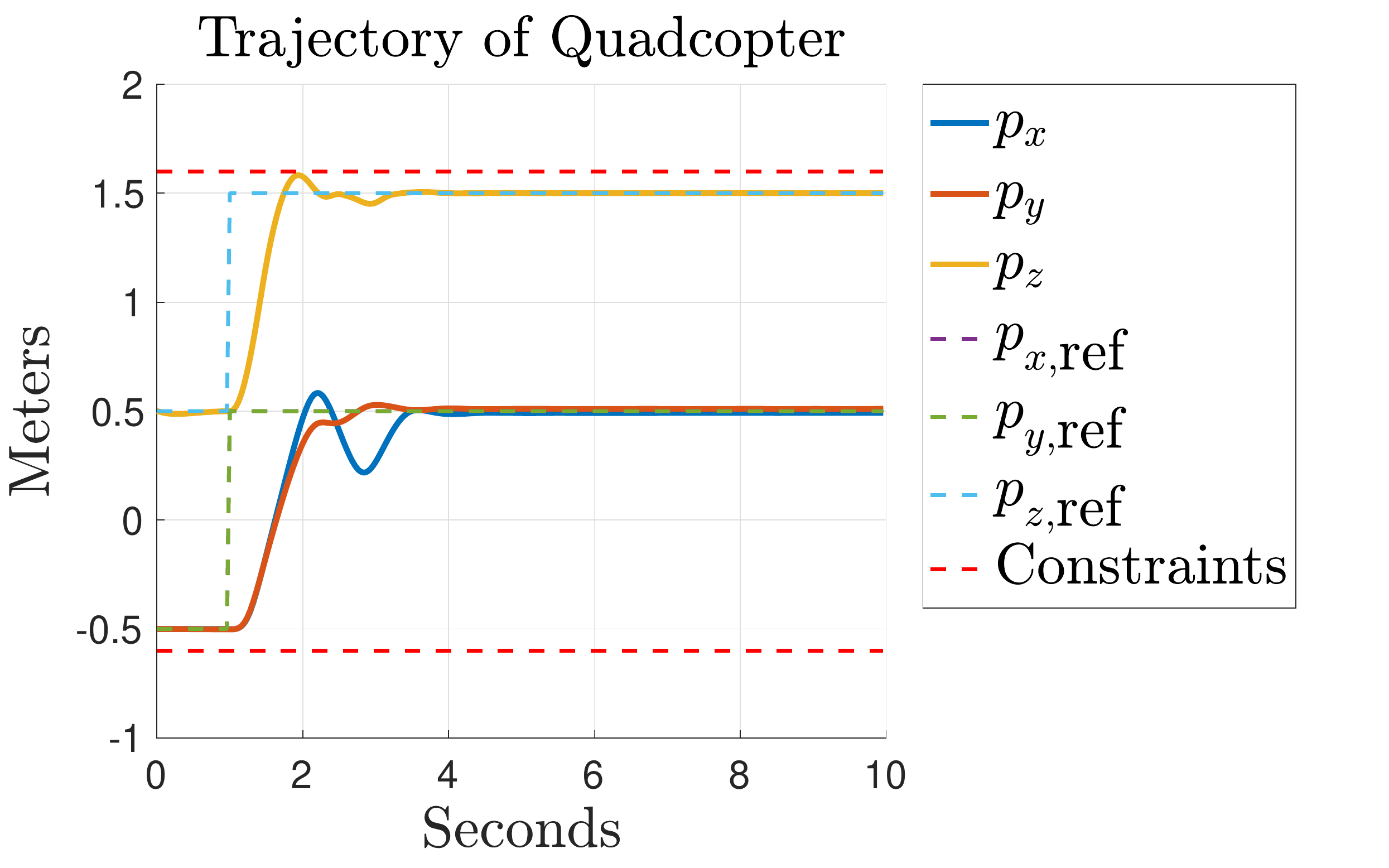}
  \caption{}
  \label{fig:trajectory}
\end{subfigure}\hfil 
\begin{subfigure}{0.33\linewidth}
  \includegraphics[width=\linewidth]{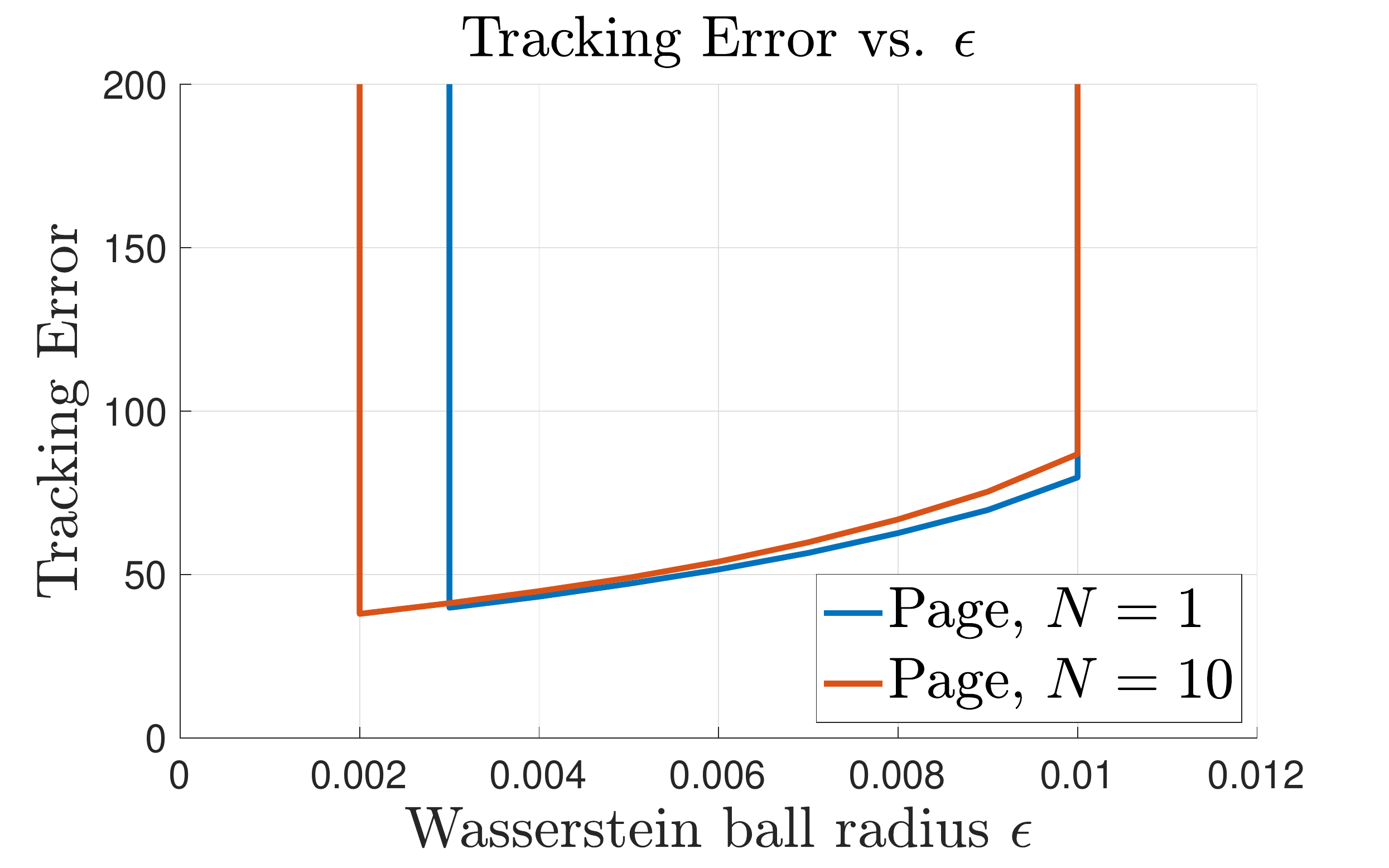}
  \caption{}
  \label{fig:epsilon_page}
\end{subfigure}\hfil
\begin{subfigure}{0.33\linewidth}
  \includegraphics[width=\linewidth]{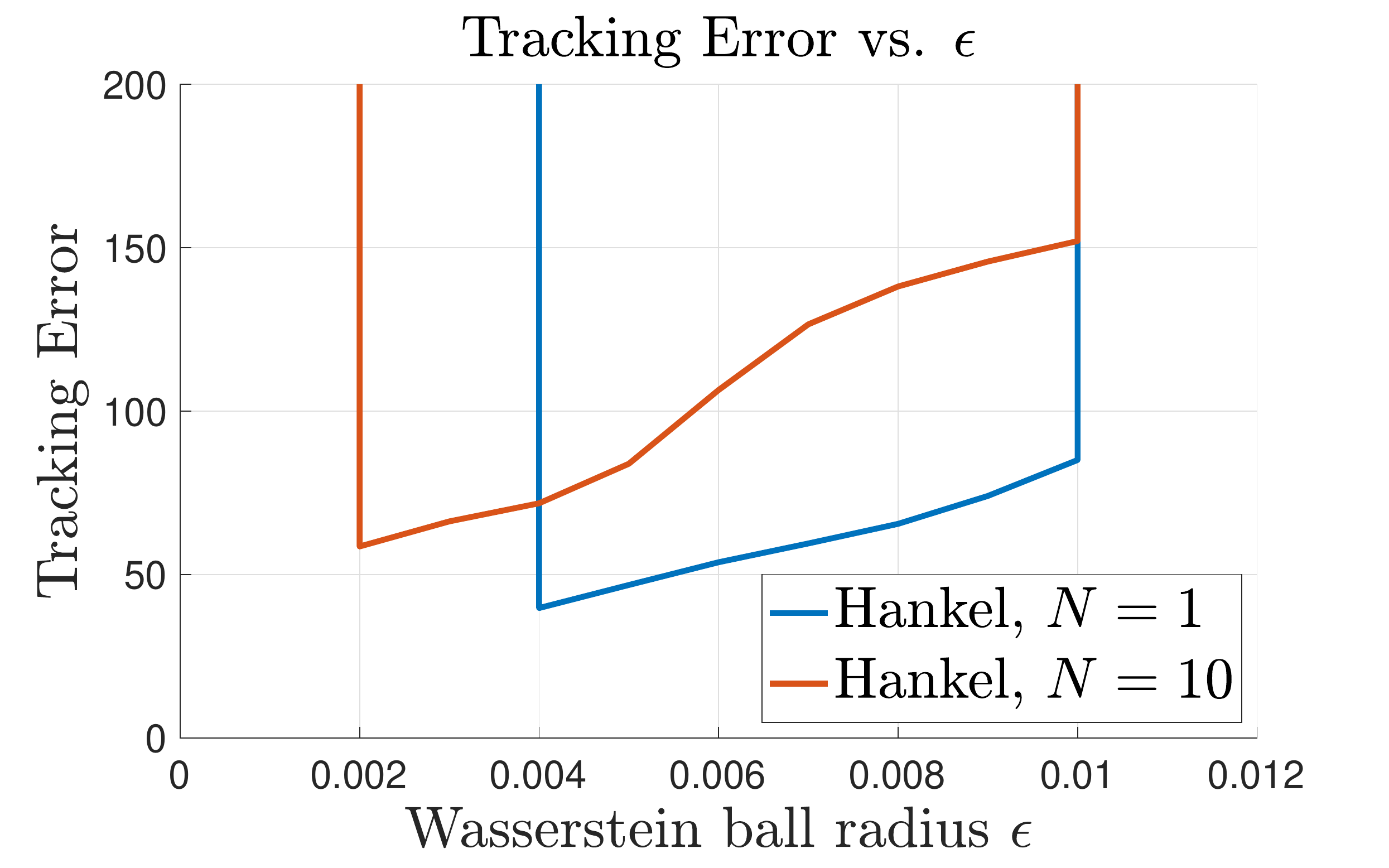}
  \caption{}
  \label{fig:epsilon_hankel}
\end{subfigure}

\medskip
\begin{subfigure}{0.33\linewidth}
  \includegraphics[width=\linewidth]{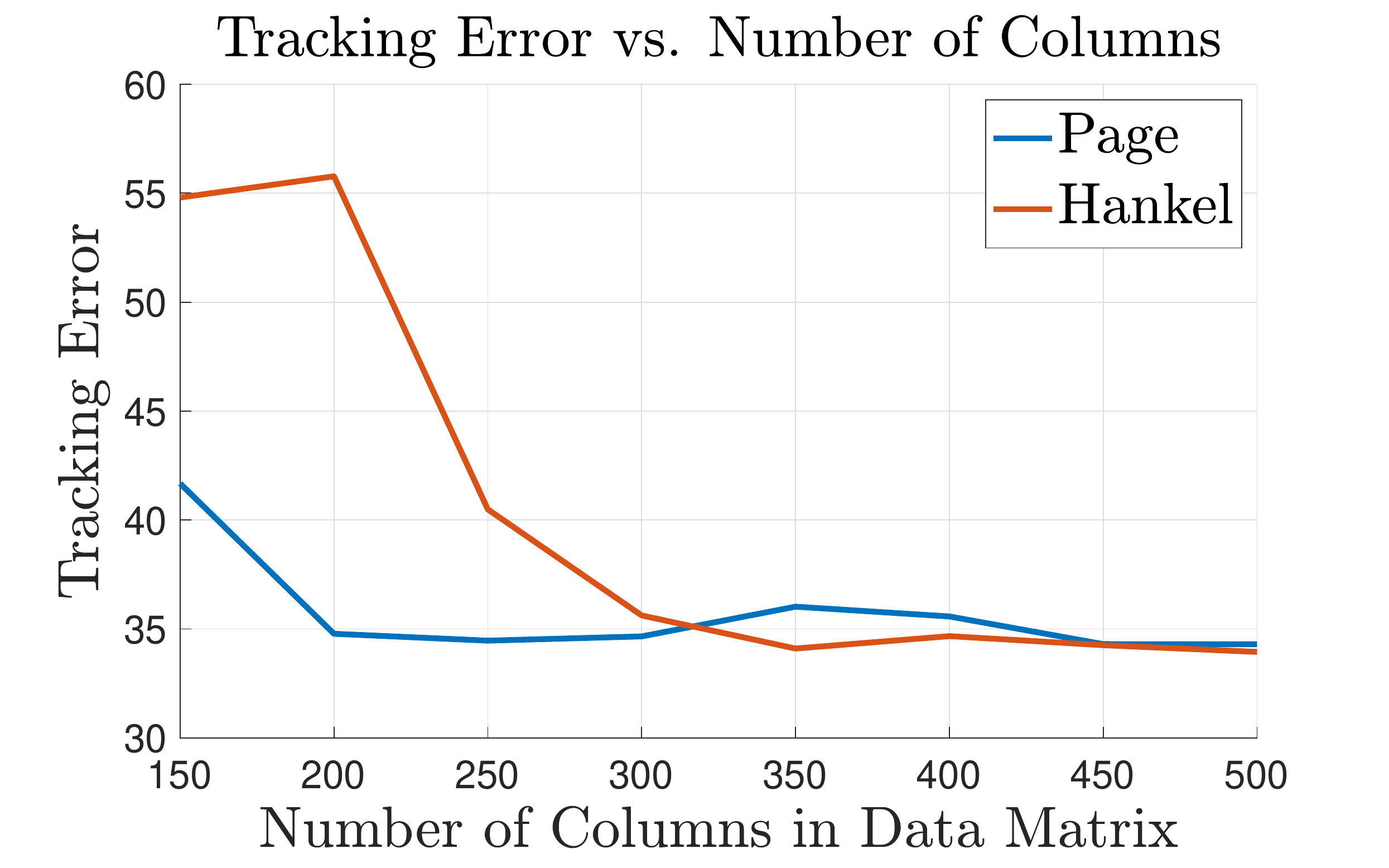}
  \caption{}
  \label{fig:number_columns}
\end{subfigure}\hfil
\begin{subfigure}{0.33\linewidth}
  \includegraphics[width=\linewidth]{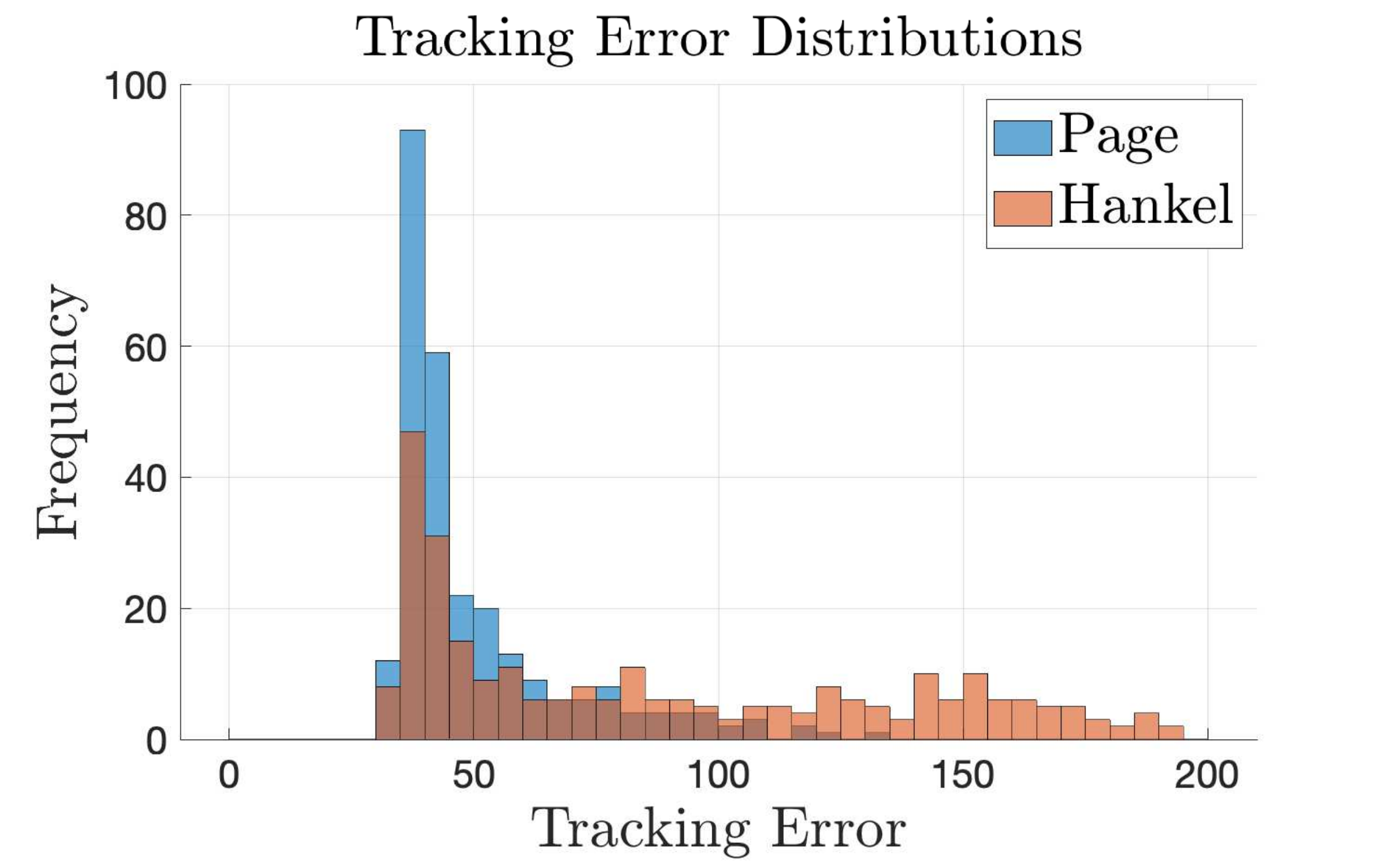}
  \caption{}
  \label{fig:page_vs_hankel}
\end{subfigure}\hfil
\begin{subfigure}{0.33\linewidth}
  \includegraphics[width=\linewidth]{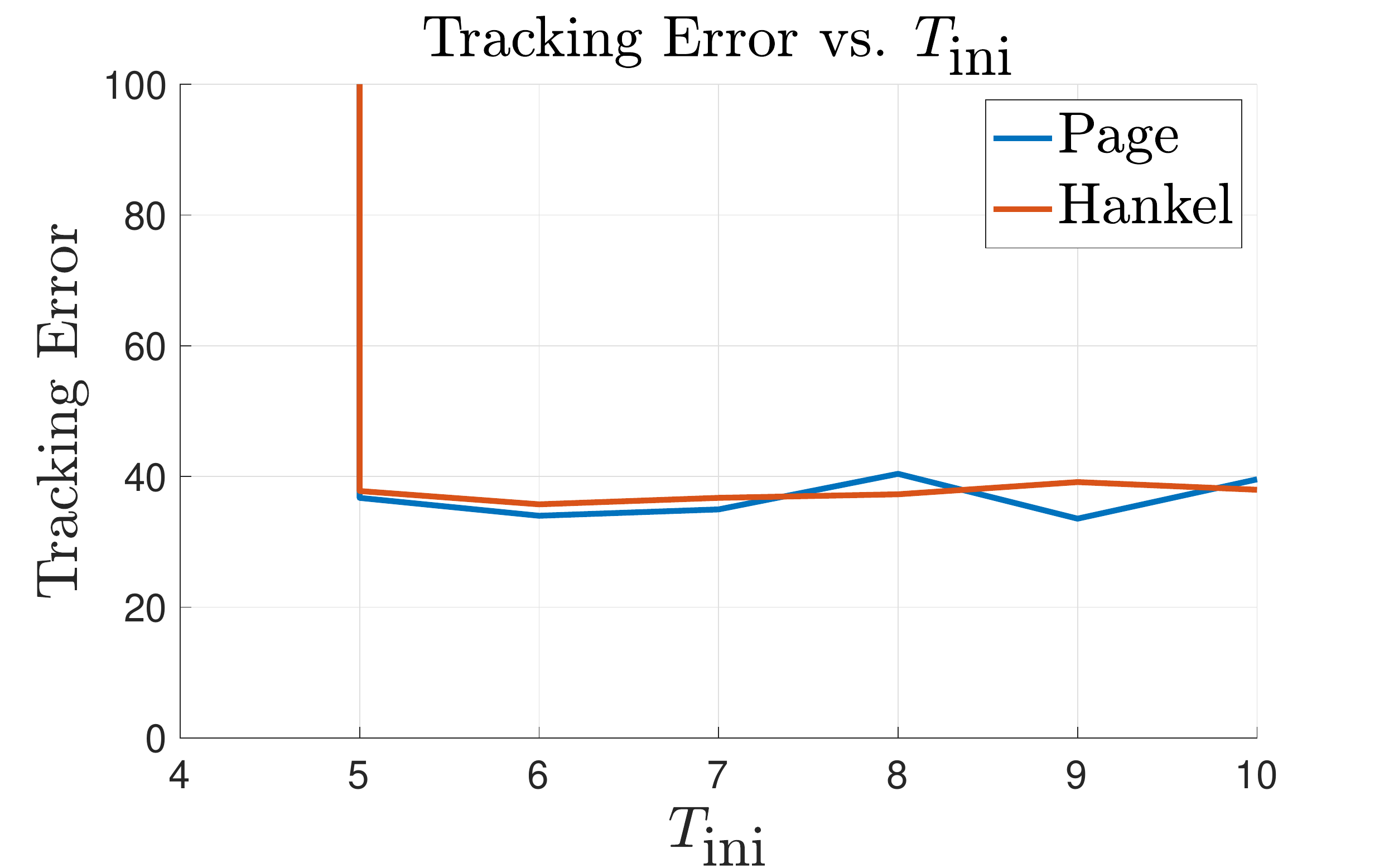}
  \caption{}
  \label{fig:Tini}
\end{subfigure}
\caption{\subref{fig:trajectory} Step trajectory of quadcopter using distributionally robust DeePC Algorithm~\ref{alg:robustdeepc}; \subref{fig:epsilon_page}--\subref{fig:epsilon_hankel} Dependence of tracking error on the Wasserstein radius $\epsilon$; \subref{fig:number_columns} Dependence of tracking error on the number of columns in the data matrix; \subref{fig:page_vs_hankel} Comparison of tracking error for Page and Hankel matrices; \subref{fig:Tini} Dependence of tracking error on the initial condition horizon $\Tini$.}
\end{figure*}

The output measurements were taken at a rate of 25Hz (i.e., every 40ms). We performed $N=10$ offline experiments each of which yielded 15500 input/output measurements to populate the Page matrices $\Uphat$, $\Ufhat$, $\Yphat^{(i)}$, and $\Yfhat^{(i)}$ for $i\in\{1,\dots,10\}$. The inputs used to excite the system were drawn from a Gaussian distribution. The same randomly generated inputs were used for every repeated experiment. Due to the fact that the quadcopter system is open-loop unstable, these randomly generated inputs were added to an existing control that maintains the quadcopter around a hover state. This existing controller was only used during the offline data collection phase, and was not used during the online control phase. Since the randomly generated inputs were much larger in magnitude than the adjustments made by the existing stabilizing controller during data collection, the matrices $\Yphat^{(i)}$, $\Yfhat^{(i)}$ are approximately i.i.d. 

The future prediction horizon was chosen as $\Tf=25$ (1 second in real time), and the time horizon used to implicitly estimate the initial condition was set to $\Tini = 6$. The inputs were constrained to the set $[0.1597,0.4791]\times [-\frac{\pi}{2}, \frac{\pi}{2}]^2$. The output constraint function $h$ was chosen such that $(p_{x},p_{y},p_{z})$ was constrained in the set $[-0.6,1.6]^3$ with $\alpha=0.1$ the confidence parameter in the CVaR constraint. The composite objective function consisted of weighted 2-norms
\[
\begin{aligned}
&f_1(u_{[1,\Tf]}) = 16\|u_{[1,\Tf],1}-f_{\textup{ref}}\|_2+4\|(u_{[1,\Tf],2},u_{[1,\Tf],3})\|_2\\
&f_2(y_{[1,\Tf]})^{2} = 1600^2{\sum\nolimits_{t=1}^{\Tf}\|y_t-y_{\textup{ref}}\|_2^2}\\
&f_3(\sigma)=750000\|\sigma\|_2,
\end{aligned}
\]
where $f_{\textup{ref}}=(0.27,\dots,0.27)$ is the thrust needed to hover the quadcopter and $y_{\textup{ref}}=(p_{x,\textup{ref}},p_{y,\textup{ref}},p_{z,\textup{ref}})=(0.5,0.5,1.5)$. The Wasserstein norm was chosen as the 2-norm and $\eps = 0.003$. The optimization problem~\eqref{eq:tractabledeepc} was implemented in a receding horizon fashion with control horizon $\nu=1$ (see Algorithm~\ref{alg:robustdeepc}). A representative trajectory of the nonlinear and stochastic closed-loop system is illustrated in Figure~\ref{fig:trajectory}, where the distributionally robust DeePC succeeds in steering the quadcopter to the reference trajectory while satisfying output constraints.

\subsection{Choosing the Hyperparameters}\label{sec:chooseparams}
In this section we study the effect of the hyperparameters on the performance of distributionally robust DeePC. To compare the performance of Algorithm~\ref{alg:robustdeepc} for different hyperparameters, we simulated the nonlinear stochastic quadcopter in receding horizon with control horizon $\nu=1$. We computed the tracking error $\sum_{t=0}^{T_{\textup{sim}}} \|y_t-y_{\textup{ref}}\|_2^2$ where $T_{\textup{sim}}=250$ (corresponding to 10 seconds in real time). 

\subsubsection{$\eps$}
We varied the Wasserstein radius for fixed $N=1$ and $N=10$, $\Tini=6$, and 300 columns in the data matrix. Figures~\ref{fig:epsilon_page} and~\ref{fig:epsilon_hankel} show the relationship between tracking error and $\eps$ for Page and Hankel matrix structures, respectively. Increasing the number of data batches $N$ increases the range of radii that lead to satisfactory tracking performance. This is because the empirical distribution $\Phat_N$ at the centre of the Wasserstein ball in problem~\eqref{eq:robustdeepc} gets closer (in the Wasserstein sense) to the true distribution $\Prob$ from which the data is drawn. This allows the user to decrease $\eps$. Since the optimal value of the distributionally robust DeePC problem~\eqref{eq:tractabledeepc} is monotonically increasing with $\eps$, choosing the minimum epsilon such that $\Prob \in\B$, gives the minimum upper bound on the out-of-sample performance that holds with high confidence. Hence, increasing $N$ decreases this minimum $\eps$ leading to better out-of-sample performance with high confidence. We note that the effect of increasing $N$ is more pronounced for the Hankel matrices. This may be due to the fact that Page matrices already require more data points to construct when compared to a Hankel matrix with the same number of columns.

Performing experiments to obtain the range of $\eps$ resulting in stable flight (Figures~\ref{fig:epsilon_page} and~\ref{fig:epsilon_hankel}) may not always be possible.
We refer the reader to~\cite[Section 7.2.2]{PME-DK:18} for methods approximating an optimal $\eps$ using the data already collected.

\subsubsection{$T$}
We varied the number of data points used in the data matrices constructed in~\eqref{eq:data}. In particular, we studied the effect of adding additional columns into the data matrices with $N=1$, $\Tini=6$, and choosing $\eps \in \{a\times 10^{-b}\mid a\in\{1,\dots,9\}, b\in\{2,3,4\}\}$ which minimizes the tracking error. In these simulations the number of data points used to construct the Page matrices did not meet the theoretical minimum required by Theorem~\ref{thm:fundamental}, but still exhibits good tracking performance. Figure~\ref{fig:number_columns} indicates that increasing the amount of data in the data matrices can significantly improve the tracking performance. Past a certain threshold (200 columns for Page matrix structure, 300 columns for Hankel matrix structure), the tracking performance remains approximately constant. This is likely due to the fact that with this amount of data, the data matrices characterize a rich enough subspace of trajectories well approximating the nonlinear dynamics. This intuition is drawn from the fact that nonlinear dynamics (under certain assumptions) can be lifted to large (often infinite) dimensional linear dynamics, and a large enough data set well approximates the dominant modes.

Aside from comparing the number of columns (and thus the size of the optimization variables), we also compare the Hankel and Page matrices given the same fixed amount of data $T$. A representative comparison is presented in Table~\ref{table:solvetime}. The mean solve time represents the average time in seconds it took to solve~\eqref{eq:tractabledeepc} with $N=1$, $\Tini=6$ and choosing $\eps \in \{a\times 10^{-b}\mid a\in\{1,\dots,9\}, b\in\{2,3,4\}\}$ which minimizes the tracking error. The optimization was performed using MOSEK~\cite{mosek} on a 3.4 GHz Intel Core i5 with 16GB of RAM.
\begin{table}[h]
\centering
\caption{Comparison of mean solve times and tracking error for Hankel and Page matrices for varying data lengths.}\label{table:solvetime}
    \begin{tabular}{| l | c | c | c | c |}    \hline
      & \multicolumn{2}{|c|}{Mean solve time (s)} & \multicolumn{2}{|c|}{Tracking error} \\ \hline
     Data length & $T=527$ & $T=15407$ & $T=527$ & $T=15407$ \\ \hline
    Hankel & 0.182 & 3.70 & 33.94 & 34.35 \\ \hline
    Page & 0.007 & 0.166 & $7\times 10^6$ & 33.97 \\ \hline
    \end{tabular}
    \end{table}

When $T=527$, the number of columns in the Hankel and Page matrices are 497 and 17 respectively. When $T=15407$, the number of columns in the Hankel and Page matrices are 15377 and 497 respectively. Each additional $\Tini+\Tf$ data points allow the addition of one column to the Page matrix and $\Tini+\Tf$ columns to the Hankel matrix. Hence, the size of the optimization variable $g$ in~\eqref{eq:tractabledeepc} when using a Hankel matrix increases $\Tini+\Tf$ times faster than that of a Page matrix.

In terms of tracking error, we see that for $T=527$ the performance of the Page matrix is poor (i.e., the quadcopter crashes) since there are not enough columns in the Page matrix to construct a rich enough subspace of trajectories. Increasing to $T=15407$ results in good performance for both the Hankel and Page matrices at the cost of a much larger solve time for the Hankel matrix. Furthermore, the tracking errors confirm the observation from Figure~\ref{fig:number_columns} that there is a threshold on the number of data points above which the tracking errors remain approximately constant.

\subsubsection{Data matrix structure}
We performed simulations with hyperparameters $N\in\{1,2,5,10\}$, $\Tini=6$, number of columns in data matrix ranging from 150 to 500, and $\eps\in \{a\times 10^{-b}\mid a\in\{1,\dots,9\}, b\in\{2,3,4\}\}$.
Over all simulations, the same data set was used to populate the Page/Hankel matrices constructed in~\eqref{eq:data}. The 273 simulations with tracking error no larger than 200 (i.e., those exhibiting stable flight of the quadcopter) are plotted in the histogram in Figure~\ref{fig:page_vs_hankel}. The Page matrix significantly outperforms the Hankel matrix in terms of tracking performance. This observation may be due to the fact that the tractable reformulation of the distributionally robust DeePC objective in~\eqref{eq:tractabledeepc} is tight for the Page matrix. Additionally, the entries of the Hankel matrix are repeated which may result in a higher sensitivity to noise compared to the Page matrix whose entries are not repeated. Lastly, we note that using a Page matrix allows one to perform singular value decomposition (SVD) to optimally preprocess the data. This is not the case for the Hankel matrix since SVD does not preserve Hankel structure. See~\cite{LH-JC-JL-FD:19-oscillation} for a relevant case study.

\subsubsection{$\Tini$}
We varied the horizon over which the initial condition is implicitly estimated while fixing $N=1$, and 300 columns in the data matrix. For each value of $\Tini\in\{1,\dots,10\}$, a simulation was performed for every $\eps\in\{a\times 10^{-b}\mid a\in\{1,\dots,9\}, b\in\{2,3,4\}\}$. 
For the best performing $\eps$, the value of the tracking error for each $\Tini$ is reported in Figure~\ref{fig:Tini}. Once $\Tini$ is past a certain threshold ($\Tini=5$), the distributionally robust DeePC algorithm exhibits satisfactory tracking performance. This is because $\Tini$ is the main parameter to fix the ``system complexity'' inside the DeePC algorithm (see Lemma~\ref{lem:initialstate}).

\section{Conclusion}\label{sec:conclusion}
We presented a data-driven method for controlling stochastic constrained LTI systems only using raw data collected from the system and without the need to explicitly identify a model. The method comes with strong out-of-sample performance guarantees and ensures constraint satisfaction with high confidence. Furthermore, we discussed how the performance of the algorithm is affected by the hyperparameters. Future work includes exploring the extension of this method to strongly nonlinear systems, extending to closed-loop guarantees, and comparison to identification-based control approaches. Furthermore, our Page excitation condition is only sufficient, and we seek relaxed conditions. Finally, we are interested in an online adaptation of the data matrix.

\section*{Acknowledgment}
The authors would like to thank Linbin Huang, Peyman Mohajerin Esfahani, Claudio de Persis, Ivan Markovsky, and Manfred Morari for useful discussions, as well as Ezzat Elokda and Paul Beuchat for providing the simulation model.

\appendix
\renewcommand{\thelemma}{A.\arabic{lemma}}
\renewcommand{\theremark}{A.\arabic{remark}}
\begin{IEEEproof}[Proof of Theorem~\ref{thm:fundamental}]
The ``if'' part of the statement is obvious by linearity and time-invariance of the system. We now prove the ``only if'' part. Our proof strategy is partially inspired by arguments in \cite{HJVW-CDP-MKC-PT:20}. Let $\bar{u}_{[1,L]} = (\bar{u}_1,\dots,\bar{u}_{L})$, $\bar{y}_{[1,L]} = (\bar{y}_1,\dots,\bar{y}_{L})$ be a trajectory of $\bv$. Hence, there exist matrices $A,B,C,D$ such that $(\bar{u}_{[1,L]},\bar{y}_{[1,L]})$ is a trajectory of $\bv(A,B,C,D)$ with initial state $\bar{x}_1\in\real^{\bm{n}(\bv)}$. For ease of notation, throughout the proof we will denote $\bm{n}(\bv)$ by $n$. Define the matrix of  state sequences from $x_{i}$ to $x_{i+jL}$ as
\[
X_{i,j}=\begin{pmatrix} x_i & x_{i+L} & \cdots & x_{i+jL}\end{pmatrix}.
\]
Then by the system dynamics we have that
\[
\begin{aligned}
&\begin{pmatrix}
\mathscr{P}_L(u_{[1,T]})\\
\mathscr{P}_L(y_{[1,T]})
\end{pmatrix}\\
&=\begin{pmatrix}
0 & I \\
\mathcal{O}_L(A,C) & \mathcal{T}_L(A,B,C,D)
\end{pmatrix}\begin{pmatrix}
X_{1,\lfloor \frac{T}{L}\rfloor-1} \\
\mathscr{P}_L(u_{[1,T]})
\end{pmatrix},
\end{aligned}
\]
where $I$ is the identity matrix. Furthermore,
\[
\begin{pmatrix}
\bar{u}_{[1,L]} \\
\bar{y}_{[1,L]}
\end{pmatrix}
=
\begin{pmatrix}
0 & I \\
\mathcal{O}_L(A,C) & \mathcal{T}_L(A,B,C,D)
\end{pmatrix}
\begin{pmatrix}
\bar{x}_1 \\
\bar{u}_{[1,L]}
\end{pmatrix}.
\]
If there exists a vector $g$ such that
\begin{equation}\label{eq:existg}
\begin{pmatrix}
\bar{x}_1 \\
\bar{u}_{[1,L]}
\end{pmatrix}
=
\begin{pmatrix}
X_{1,\lfloor \frac{T}{L}\rfloor-1} \\
\mathscr{P}_L(u_{[1,T]})
\end{pmatrix}g,
\end{equation}
then
\[
\begin{aligned}
\begin{pmatrix}
\bar{u}_{[1,L]} \\
\bar{y}_{[1,L]}
\end{pmatrix}
&=
\begin{pmatrix}
0 & I \\
\mathcal{O}_L(A,C) & \mathcal{T}_L(A,B,C,D)
\end{pmatrix}
\begin{pmatrix}
X_{1,\lfloor \frac{T}{L}\rfloor-1} \\
\mathscr{P}_L(u_{[1,T]})
\end{pmatrix}g\\
&=
\begin{pmatrix}
\mathscr{P}_L(u_{[1,T]})\\
\mathscr{P}_L(y_{[1,T]})
\end{pmatrix}g,
\end{aligned}
\]
which is what we wanted to show. Hence it remains to prove that there always exists a vector $g$ satisfying~\eqref{eq:existg} for any $\bar{x}_1$ and $\bar{u}_{[1,L]}$. Equivalently, we will show that the matrix $\begin{pmatrix}
X_{1,\lfloor \frac{T}{L}\rfloor-1} \\
\mathscr{P}_L(u_{[1,T]})
\end{pmatrix}$ has full row rank.
Let $\begin{pmatrix}\xi & \eta\end{pmatrix}$ be an arbitrary vector in the left kernel of this matrix where $\xi^\top \in \real^n$ and $\eta^\top\in\real^{mL}$. We will show that $\begin{pmatrix}\xi & \eta\end{pmatrix}$ must be the zero vector. To do this, we switch our attention to the matrix
\begin{equation}\label{eq:XU}
\begin{pmatrix}
X \\
U
\end{pmatrix}\defeq \begin{pmatrix}
X_{1,\lfloor \frac{T}{L}\rfloor-n-1} \\
\mathscr{P}_L(u_{[1,T-nL]}) \\
\mathscr{P}_L(u_{[L+1,T-(n-1)L]}) \\
\vdots \\
\mathscr{P}_L(u_{[nL+1,T]}) 
\end{pmatrix}
\end{equation}
and follow arguments from \cite{HJVW-CDP-MKC-PT:20} to construct $n+1$ vectors in the left kernel of~\eqref{eq:XU}. By definition, $\begin{pmatrix} \xi & \eta & \vectorzeros[nmL]\end{pmatrix}$ is in the left kernel of~\eqref{eq:XU}, where $\vectorzeros[nmL]^\top\in\real^{nmL}$ denotes a vector containing $nmL$ zeros.
\begin{figure*}[b]
\hrulefill
\begin{equation}
M\defeq\begin{pmatrix}
A^L & A^{L-1}B & \cdots & AB & B & \vectorzeros[{n\times mL}] & \cdots & \vectorzeros[{n\times mL}] \\
\vectorzeros[{m\times n}] & \vectorzeros[{m\times m}] & \cdots & \vectorzeros[{m\times m}] & \vectorzeros[{m\times m}] & \begin{bmatrix} I & \cdots & \vectorzeros[{m\times m}]\end{bmatrix} & \cdots & \begin{bmatrix} \vectorzeros[{m\times m}] & \cdots & \vectorzeros[{m\times m}]\end{bmatrix}\\
\vdots & \vdots & \ddots & \vdots & \vdots & \vdots & \ddots & \vdots \\
\vectorzeros[{m\times n}] & \vectorzeros[{m\times m}] & \cdots & \vectorzeros[{m\times m}] & \vectorzeros[{m\times m}] & \begin{bmatrix} \vectorzeros[{m\times m}] & \cdots & I \end{bmatrix} & \cdots & \begin{bmatrix} \vectorzeros[{m\times m}] & \cdots & \vectorzeros[{m\times m}]\end{bmatrix}
\end{pmatrix}
\label{eq:M}
\end{equation}
\end{figure*}
Furthermore, we have that
\begin{equation}
\begin{pmatrix}
X_{L+1,\lfloor \frac{T}{L}\rfloor-n-1} \\
\mathscr{P}_L(u_{[L+1,T-(n-1)L]})
\end{pmatrix}=
M
\begin{pmatrix}
X \\
U
\end{pmatrix},
\label{eq:aux-equation}
\end{equation}
where $M$ is defined in~\eqref{eq:M} and $\vectorzeros[n\times n]$ is the zero matrix.
Since $\begin{pmatrix}
X_{L+1,\lfloor \frac{T}{L}\rfloor-n-1} \\
\mathscr{P}_L(u_{[L+1,T-(n-1)L]})
\end{pmatrix}$ is a submatrix of $\begin{pmatrix}
X_{1,\lfloor \frac{T}{L}\rfloor-1} \\
\mathscr{P}_L(u_{[1,T]})
\end{pmatrix}$ then by definition of $\begin{pmatrix}\xi & \eta \end{pmatrix}$, we have
\[
\begin{pmatrix} \xi & \eta \end{pmatrix}\begin{pmatrix}
X_{L+1,\lfloor \frac{T}{L}\rfloor-n-1} \\
\mathscr{P}_L(u_{[L+1,T-(n-1)L]})
\end{pmatrix}= \vectorzeros[(\lfloor \frac{T}{L}\rfloor -n-1)L+1].
\]
Hence, by \eqref{eq:aux-equation} we have that
\[
\begin{pmatrix} \xi A^L & \xi A^{L-1}B & \cdots & \xi AB & \xi B & \eta & \vectorzeros[{(n-1)mL}]\end{pmatrix}
\]
is in the left kernel of~\eqref{eq:XU}. Note also that 
\[
\begin{aligned}
&X_{2L+1,\lfloor \frac{T}{L}\rfloor-n-1} \\
&= \begin{pmatrix} A^{2L} & A^{2L-1}B & \cdots & AB & B & \vectorzeros[{n\times (n-1)mL}]\end{pmatrix}
\begin{pmatrix}
X \\
U
\end{pmatrix}.
\end{aligned}
\]
Proceeding as above we find that the vector
\setcounter{MaxMatrixCols}{20}
\[
\begin{pmatrix} 
\xi A^{2L} & \xi A^{2L-1}B & \cdots & \xi AB & \xi B & \eta & \vectorzeros[{(n-2)mL}]
\end{pmatrix}
\]
is in the left kernel of~\eqref{eq:XU}. Continuing in this way, we obtain $n+1$ vectors in the left kernel of~\eqref{eq:XU} given by
\[
\begin{aligned}
&w_0 = \begin{pmatrix} \xi & \eta & \vectorzeros[{nmL}]\end{pmatrix} \\
&w_1 = \begin{pmatrix} \xi A^L & \xi A^{L-1}B & \cdots & \xi AB & \xi B & \eta & \vectorzeros[{(n-1)mL}]\end{pmatrix} \\
&w_2 = \begin{pmatrix} 
\xi A^{2L} & \xi A^{2L-1}B & \cdots & \xi AB & \xi B & \eta & \vectorzeros[{(n-2)mL}]
\end{pmatrix} \\
&\vdots \\
&w_n =\begin{pmatrix} 
\xi A^{nL} & \xi A^{nL-1}B & \cdots & \xi AB & \xi B & \eta
\end{pmatrix}
\end{aligned}
\]
However, since the input vectors are $L$-Page exciting of order $n+1$ the left kernel of~\eqref{eq:XU} has dimension at most $n$. This implies that the vectors $w_0,\dots, w_n$ are linearly dependent. From the structure of the vectors we see that $\eta=\vectorzeros[{mL}]$. By Cayley-Hamilton, there exist $\alpha_i\in\real$, $i\in\{0,\dots,n\}$ with $\alpha_n=1$ such that $\sum_{i=0}^n \alpha_i (A^L)^i=0$.
Define
\[
\begin{aligned}
&v = \sum_{i=0}^n \alpha_i w_i\\
&=\begin{pmatrix}
\xi \sum\limits_{i=0}^n\alpha_i (A^L)^i \!\!\!& \cdots &\!\!\! \xi  A^{L-1} B \alpha_n \!\!\!& \cdots &\!\!\! \xi B \alpha_n & \vectorzeros[mL] \end{pmatrix}\\
&=\begin{pmatrix}
\vectorzeros[n] & \xi \sum\limits_{i=1}^n\alpha_i  A^{iL-1} B \!\!\!& \cdots &\!\!\! \xi  A^{L-1} B \!\!\!& \cdots &\!\!\! \xi B & \vectorzeros[mL] \end{pmatrix}.
\end{aligned}
\]
Since $v$ is a linear combination of vectors $w_i$ then it is itself in the left kernel of~\eqref{eq:XU}. Hence, $\begin{pmatrix} \xi \sum_{i=1}^n\alpha_i  A^{iL-1} B & \cdots & \xi  A^{L-1} B & \cdots & \xi B & \vectorzeros[mL] \end{pmatrix}$ is in the left kernel of $U$. However, $U$ has full row rank by the excitation assumption implying the above vector must be the zero vector. Thus, $\xi \begin{pmatrix} A^{L-1}B & \cdots & B\end{pmatrix}=\vectorzeros[mL]$. By assumption we have that $L\geq n$ and $\bv$ is controllable. Thus $\xi=\vectorzeros[n]$. Hence, $\begin{pmatrix} \xi & \eta\end{pmatrix}=\vectorzeros[n+mL]$ implying that~\eqref{eq:XU} has full row rank which concludes the proof.
\end{IEEEproof}
\begin{IEEEproof}[Proof of Proposition~\ref{thm:setequiv}]
We first look at the feasible set of~\eqref{eq:mpc}. By rewriting the constraints in~\eqref{eq:mpc} we obtain
\[
y= \mathcal{O}_{\Tf}(A,C)\hat{x}_t+\mathcal{T}_{\Tf}(A,B,C,D) u,
\quad u\in\mathcal{U},\; y\in\mathcal{Y},
\]
where $\hat{x}_t$ is the state at time $t$. We now look at the feasible set of~\eqref{eq:deepc}. Since $\hat{u}_{[1,T]}$ is $(\Tini+\Tf)$-Page exciting of order $\bm{n}(\bv)+1$, then by Theorem~\ref{thm:fundamental} $\textup{image}\left(\col(\Uphat,\Yphat,\Ufhat,\Yfhat)\right)=\bv_{\Tini+\Tf}$, where $\Uphat$, $\Yphat$, $\Ufhat$, $\Yfhat$ are defined as in~\eqref{eq:data}. Hence, the feasible set of~\eqref{eq:deepc} is equal to $\{(u,y)\in \mathcal{U}\times\mathcal{Y}\mid \col(\uinihat,u,\yinihat,y)\in\bv_{\Tini+\Tf}\}$. Since the system $\bv$ yields an equivalent representation given by $\bv(A,B,C,D)$, then by Lemma~\ref{lem:initialstate} the feasible set of~\eqref{eq:deepc} can be written as the set of pairs $(u,y)\in \mathcal{U}\times\mathcal{Y}$ satisfying
\[
y= \mathcal{O}_{\Tf}(A,C)\xinihat+\mathcal{T}_{\Tf}(A,B,C,D) u,
\]
where $\xinihat$ is uniquely determined from $\col(\uinihat,\yinihat)$ and hence coincides with $\hat{x}_t$. This proves the claim.
\end{IEEEproof}

The following two lemmas will serve as essential cornerstones of our subsequent results.
\begin{lemma}\label{lem:dual}
Assume that $\varphi$ is proper, convex, and lower semicontinuous. Then
\[
\sup_{\Q\in\B}\Exp_\Q[\varphi(\xi)]=\inf_{\lambda\geq 0}\lambda\eps+\frac{1}{N} \sum_{i=1}^N \sup_{\xi\in\Xi}(\varphi(\xi)-\lambda\|\xi-\hat{\xi}^{(i)}\|_r)
\]
\end{lemma}
\begin{IEEEproof}
The proof follows from~\cite[Theorem 4.2]{PME-DK:18}, specifically from Equation (12b) in~\cite{{PME-DK:18}} which relies on a marginalization and dualization of the constraint $\Q\in\B$.
\end{IEEEproof}
Next we present a similar result as~\cite[Lemma 4.1]{JC-JL-FD:19}. The result relaxes the conservative assumption that $\Xi=\prod_{i=1}^{p(\Tini+\Tf)}\Xi_i$ for appropriately defined $\Xi_i$ imposed in~\cite{JC-JL-FD:19}. The proof strategy is partially inspired by arguments in~\cite{PME-DK:18}.
\begin{lemma}\label{lem:lipschitz}
Let $c,d \in \integerspositive$ and $\Omega \subseteq \real^{cd}$. Let $q$ be such that $\frac{1}{r}+\frac{1}{q}=1$. Let $\hat{\zeta} = \col(\hat{\zeta}_1,\dots,\hat{\zeta}_c)\in\Omega$ be given where each $\hat{\zeta}_i\in\real^d$. Let $\lambda\in\realpositive$ and $\varphi\colon \real^c\to \real$ be a convex and Lipschitz continuous function with Lipschitz constant $L_{\varphi}$ with respect to the $r$-norm. Then for $b\in\real^d$,
\[
\begin{aligned}
&\sup_{\zeta \in \Omega} \varphi(\zeta_1^{\top}b,\dots, \zeta_c^{\top}b) -\lambda\|\zeta-\hat{\zeta}\|_r\\
&\leq \begin{cases}
\varphi(\hat{\zeta}_1^{\top}b,\dots, \hat{\zeta}_c^{\top}b) & \textup{if }L_{\varphi}\|b\|_q\leq \lambda \\
\infty & \textup{otherwise}
\end{cases}
\end{aligned}
\]
The above is an equality when $\Omega=\real^{cd}$.
\end{lemma}
\begin{IEEEproof}
We begin by first noting that 
\[
\begin{aligned}
&\sup_{\zeta \in \Omega} \varphi(\zeta_1^{\top}b,\dots, \zeta_c^{\top}b) -\lambda\|\zeta-\hat{\zeta}\|_r\\
&\leq \sup_{\zeta \in \real^{cd}} \varphi(\zeta_1^{\top}b,\dots, \zeta_c^{\top}b) -\lambda\|\zeta-\hat{\zeta}\|_r 
\end{aligned}
\]
with equality when $\Omega =\real^{cd}$. Let $\Phi(\zeta)=\varphi(\zeta_1^{\top}b,\dots,\zeta_c^{\top}b)$.
By definition of the conjugate function,
\begin{align*}
\Phi^{\ast}(z)&=\underset{\zeta\in \real^{cd}}\sup \langle z,\zeta \rangle- \Phi(\zeta)\\
&=\underset{\zeta\in \real^{cd}}\sup \langle z,\zeta \rangle- \varphi(\zeta_1^{\top}b,\dots,\zeta_c^{\top}b)\\
 &=\underset{\zeta\in \real^{cd},s\in\real^c}\sup\left\{\langle z,\zeta\rangle- \varphi(s)\;\middle\vert \;\underset{\forall j\in\{1,\dots,c\}}{s_j=\zeta_j^{\top} b}\right\},
\end{align*}
where $s=\col(s_1,\dots,s_c)$.
The Lagrangian of the above is given by
\[
\mathscr{L}(s,\zeta,\theta)=\langle z,\zeta\rangle+\langle \theta,(s_1-\zeta_1^{\top} b,\dots, s_c-\zeta_c^{\top} b)\rangle - \varphi(s),
\]
By strong duality (see, e.g.,~\cite[Proposition 5.3.1]{DPB:09-book}),
\[
\begin{aligned}
 \Phi^{\ast}(z)&=\inf_{\theta}\sup_{\zeta\in\real^{cd},s}\langle z,\zeta\rangle \\
 &\quad\quad+\langle \theta,(s_1-\zeta_1^{\top} b,\dots,s_c-\zeta_c^{\top} b)\rangle - \varphi(s)\\
&=\inf_{\theta}\sup_{\zeta\in\real^{cd}} \langle z,\zeta\rangle -\langle \theta,(\zeta_1^{\top}b,\dots,\zeta_c^{\top}b)\rangle + \varphi^{\ast}(\theta)\\
&=\inf_{\theta}\sup_{\zeta\in\real^{cd}} \langle (z_1-\theta_1b,\dots,z_c-\theta_cb),\zeta\rangle +\varphi^{\ast}(\theta)
\end{aligned}
\]
where $z=\col(z_1,\dots,z_c)$. 
The above problem may take the value $\infty$ unless $z_i=\theta_i b$ for all $i\in\{1,\dots,c\}$. Hence, by taking $\Theta\defeq \{\theta\;\vert\; \varphi^{\ast}(\theta)<\infty\}$ as the effective domain of the conjugate function of $\varphi$, we obtain
\begin{align*}
\Phi^{\ast}(z)&= \inf_{\theta}\left\{\varphi^{\ast}(\theta)\;\middle\vert\; z_i=\theta_i b,\; \forall i\in\{1,\dots,c\}\right\}\\
&=\inf_{\theta\in\Theta}\left\{\varphi^{\ast}(\theta)\;\middle\vert\; z_i=\theta_i b,\; \forall i\in\{1,\dots,c\}\right\}.
\end{align*}
Since $\Phi$ is convex and continuous, the biconjugate $\Phi^{\ast\ast}$ coincides with the function $\Phi$ itself. Hence,
\begin{align*}
\Phi(\zeta)&=\sup_{z} \langle z,\zeta\rangle-\Phi^{\ast}(z)\\
&= \sup_{z}\left\{ \langle z,\zeta\rangle-\inf_{\theta\in\Theta}\varphi^{\ast}(\theta)\middle\vert z_i=\theta_i b,\; \forall i\in\{1,\dots,c\}\right\}\\
&=\sup_{\theta\in\Theta} \langle (\theta_1b,\dots,\theta_c b),\zeta\rangle -\varphi^{\ast}(\theta).
\end{align*}
Hence,
\begin{align*}
&\sup_{\zeta\in\real^{cd}}\Phi(\zeta)-\lambda\|\zeta-\hat{\zeta}\|_r\\
&=\sup_{\zeta\in\real^{cd}}\sup_{\theta\in\Theta} \langle (\theta_1b,\dots,\theta_c b),\zeta\rangle -\varphi^{\ast}(\theta)-\lambda\|\zeta-\hat{\zeta}\|_r\\
&=\sup_{\zeta\in\real^{cd}}\sup_{\theta\in\Theta}\inf_{\|\mu\|_q\leq\lambda}\langle (\theta_1b,\dots,\theta_c b),\zeta\rangle -\varphi^{\ast}(\theta)-\langle \mu,\zeta-\hat{\zeta}\rangle,
\end{align*}
where the last equality comes from the definition of the dual norm and homogeneity of the norm. Using the minimax theorem (see, e.g.,~\cite[Proposition 5.5.4]{DPB:09-book}) we switch the supremum and infimum in the above giving
\begin{align*}
&\sup_{\zeta\in\real^{cd}}\Phi(\zeta)-\lambda\|\zeta-\hat{\zeta}\|_r\\
&= \sup_{\theta\in\Theta}\inf_{\|\mu\|_q\leq\lambda}\sup_{\zeta\in\real^{cd}}\langle (\theta_1b,\dots,\theta_c b),\zeta\rangle -\varphi^{\ast}(\theta)-\langle \mu,\zeta-\hat{\zeta}\rangle\\
&=\sup_{\theta\in\Theta}\inf_{\|\mu\|_q\leq\lambda}\sup_{\zeta\in\real^{cd}}\langle (\theta_1b-\mu_1,\dots,\theta_c b-\mu_c),\zeta\rangle \\
&\quad\quad\quad\quad\quad\quad\quad\quad\quad\quad-\varphi^{\ast}(\theta)+\langle \mu,\hat{\zeta}\rangle,
\end{align*}
where $\mu=\col(\mu_1,\dots,\mu_c)$. Carrying out the supremum over $\zeta$ yields
\begin{align*}
&\sup_{\zeta\in\real^{cd}}\Phi(\zeta)-\lambda\|\zeta-\hat{\zeta}\|_r\\
&= \sup_{\theta\in\Theta}\inf_{\|\mu\|_q\leq\lambda}\begin{cases}
\langle (\theta_1b,\dots,\theta_c b),\hat{\zeta}\rangle -\varphi^{\ast}(\theta) & \!\!\!\!\textup{if } \mu_i=\theta_i b \\ 
& \!\!\!\!\forall i\in\{1,\dots,c\}\\
\infty &\!\!\!\!\textup{otherwise}
\end{cases}\\
&= \sup_{\theta\in\Theta}\begin{cases}
\!\langle (\theta_1b,\dots,\theta_c b),\hat{\zeta}\rangle -\varphi^{\ast}(\theta) & \!\!\!\!\textup{if } \|(\theta_1 b,\dots,\theta_c b)\|_q\leq \lambda\\
\infty &\!\!\!\!\textup{otherwise}
\end{cases}\\
&=\begin{cases}
\varphi(\hat{\zeta}_1^{\top} b,\dots,\hat{\zeta}_c^{\top} b) & \textup{if } \sup_{\theta\in\Theta}\|(\theta_1 b,\dots,\theta_c b)\|_q\leq \lambda\\
\infty &\textup{otherwise}
\end{cases}\,,
\end{align*}
where we used the definition of the biconjugate function and $\varphi^{\ast\ast} = \varphi$ (due to convexity and continuity).
Note that $\|(\theta_1 b,\dots,\theta_c b)\|_q=\|\theta\|_q\|b\|_q$. 
From~\cite[Proposition 6.5]{PME-DK:18}, we know that $\|\theta\|_q\leq L_{\varphi}$. In fact, one can show that $\sup_{\theta\in\Theta}\|\theta\|_q= L_{\varphi}$~\cite[Remark 3]{DK-PME-VN-SS:19}. 
Substituting this into the expression above proves the result.
\end{IEEEproof}
The following lemmas enable our main results and provide separate reformulations for the objective function and constraints. The proofs strategies are partially inspired by~\cite{PME-DK:18,SSA-DK-PME:19,AH-AC-JL:18}.
\begin{lemma}\longthmtitle{Objective Reformulation}\label{lem:costreformulation}
Assume that $f_2$ and $f_3$ are convex and Lipschitz continuous functions. Specifically, let $L_\text{obj}>0$ be the Lipschitz constant with respect to the $r$-norm of the  mapping $(x,y)\mapsto f_2(x)+f_3(y)$. Then
\[
\begin{aligned}
&\underset{\Q\in\B}{\sup} \Exp_{\Q}\left[f_1(\Ufhat g) + f_2(\Yf g) + f_3(\Yp g -\yinihat) \right] \leq \\
&f_1(\Ufhat g)  \!+\!\frac{1}{N}\sum_{i=1}^N \left(f_2(\Yfhat^{(i)}g)\!+\!f_3(\Yphat^{(i)} g-\yinihat)\right) \!+\! L_{\textup{obj}}\eps \|g\|_{r,\ast}
\end{aligned}
\]
The above is an equality when $\Xi=\real^{p(\Tini+\Tf){\lfloor \frac{T}{\Tini+\Tf}\rfloor}}$.
\end{lemma}
\smallskip
\begin{IEEEproof}
By definition of $\xi$, we have $f_2(\Yf g) + f_3(\Yp g -\yinihat)=f_2(\xi_{p\Tini+1}^{\top}g,\dots,\xi_{p(\Tini+\Tf)}^{\top}g)+f_3((\xi_1^{\top}g,\dots,\xi_{p\Tini}^{\top}g)-\yinihat)$. Define
\[
\begin{aligned}
\varphi(\xi_1^{\top}g,\dots,\xi_{p(\Tini+\Tf)}^{\top}g)\defeq &f_2(\xi_{p\Tini+1}^{\top}g,\dots,\xi_{p(\Tini+\Tf)}^{\top}g) \\
& \quad+f_3((\xi_1^{\top}g,\dots,\xi_{p\Tini}^{\top}g)-\yinihat)
\end{aligned}
\]
By definition, we have $L_{\varphi}=L_{\textup{obj}}$. Define also $\Phi(\xi)\defeq\varphi(\xi_1^{\top}g,\dots,\xi_{p(\Tini+\Tf)}^{\top}g)$. Then from Lemma~\ref{lem:dual} applied to $\Phi$ we obtain
\[
\begin{aligned}
&\underset{\Q\in\B}{\sup} \Exp_{\Q}\left[f_1(\Ufhat g) + f_2(\Yf g) + f_3(\Yp g -\yinihat) \right]\\
&=\inf_{\lambda\geq 0}\lambda\eps+f_1(\Ufhat g)\\
&\quad+\frac{1}{N} \sum_{i=1}^N \sup_{\xi\in\Xi}(\varphi(\xi_1^{\top}g,\dots,\xi_{p(\Tini+\Tf)}^{\top}g)-\lambda\|\xi-\hat{\xi}^{(i)}\|_r)
\end{aligned}
\]
Applying Lemma~\ref{lem:lipschitz} yields
\[
\begin{aligned}
&\underset{\Q\in\B}{\sup} \Exp_{\Q}\left[f_1(\Ufhat g) + f_2(\Yf g) + f_3(\Yp g -\yinihat) \right]\leq\\
& \eps L_{\varphi}\|g\|_{r,\ast}\!+ \!f_1(\Ufhat g) \!+ \!\frac{1}{N} \sum_{i=1}^N \left( f_2(\Yfhat^{(i)}g) \!+\! f_3(\Yphat^{(i)} g-\yinihat)\right).
\end{aligned}
\]
The above is an equality when $\Xi=\real^{p(\Tini+\Tf)\lfloor \frac{T}{\Tini+\Tf}\rfloor}$. Substituting $L_{\varphi}=L_{\textup{obj}}$ yields the result.
\end{IEEEproof}
The following two results depend on the constraint set
\[
G_{\textup{CVaR}}\defeq\left\{g\;\middle\vert\;\sup_{\Q\in\B}\textup{CVaR}^{\Q}_{1-\alpha}\left(h(\Yf g)\right)\leq0\right\}.
\]
\begin{lemma}\longthmtitle{Lipschitz Constraint Reformulation}\label{lem:constraintreformulation_lipschitz}
Let $h$ be convex and Lipschitz continuous with Lipschitz constant $L_\text{con}$ with respect to the $r$-norm. Then
\[
\left\{ g\;\middle\vert
\begin{array}{l}
\exists \tau,\lambda,s_i\; \textup{such that}\\
-\tau \alpha + L_\text{con}\eps\|g\|_{r,\ast} +\frac{1}{N}\sum_{i=1}^N s_i\leq 0\\
\tau + h(\Yfhat^{(i)} g) \leq s_i \\
s_i\geq 0 \\
\forall i\leq N
\end{array}
\right\} \subseteq G_{\textup{CVaR}} \,.
\]
The sets coincide if $\Xi=\real^{p(\Tini+\Tf){\lfloor \frac{T}{\Tini+\Tf}\rfloor}}$ and $h(\Yf g)$ is bounded on $\Xi$ for all $g$.
\end{lemma}
\smallskip
\begin{remark}\label{rem:tight}
Note that if $\Xi=\real^{p(\Tini+\Tf){\lfloor \frac{T}{\Tini+\Tf}\rfloor}}$, then $h(\Yf g)$ can only be bounded on $\Xi$ for all $g$ if $h$ is constant. Thus, for non-constant constraint functions $h$, the above will be an inner approximation of $G_{\textup{CVaR}}$.\oprocend
\end{remark}

\begin{IEEEproof}[Proof of Lemma~\ref{lem:constraintreformulation_lipschitz}]
By definition of CVaR,
\[
\textup{CVaR}_{1-\alpha}^{\Prob}(h(y))\leq 0 \iff \inf_{\tau\in\real}-\tau\alpha+\Exp_{\Prob}[(h(y)+\tau)_{+}] \leq 0
\]
This can be shown by multiplying the definition of CVaR by $\alpha>0$ and changing $\tau$ to $-\tau$.
Letting $\ell(g,\xi)\defeq h(\xi_{p\Tini+1}^{\top} g,\dots,\xi_{p(\Tini+\Tf)}^{\top} g)$, the constraint $g\in G_{\textup{CVaR}}$ is equivalent to $g$ satisfying
\begin{equation}\label{eq:cvar}
\sup_{\Q\in\B}\inf_{\tau\in\real}-\tau\alpha+\Exp_{\Q}\left[\left(\ell(g,\xi)+\tau\right)_+\right]\leq 0.
\end{equation}
We have
\begin{equation}\label{eq:minmax}
\begin{aligned}
\sup_{\Q\in\B}\inf_{\tau\in\real}-\tau\alpha+\Exp_{\Q}[( \ell(g,\xi)+\tau)_+]\\
\leq \inf_{\tau\in\real} -\tau\alpha+\sup_{\Q\in\B}\Exp_{\Q}[( \ell(g,\xi)+\tau)_+]
\end{aligned}
\end{equation}
by the max-min inequality. From Lemma~\ref{lem:dual} we have
\[
\begin{aligned}
&\underset{\Q\in\B}{\sup}\Exp_{\Q}[(\ell(g,\xi)+\tau)_+]\\
&= \inf_{\lambda\geq 0}\lambda\eps+\frac{1}{N} \sum_{i=1}^N \sup_{\xi\in\Xi}((\ell(g,\xi)+\tau)_+-\lambda\|\xi-\hat{\xi}^{(i)}\|_r).
\end{aligned}
\]
Recalling that $\ell(g,\xi)=h(\xi_{p\Tini+1}^{\top}g,\dots,\xi_{p(\Tini+\Tf)}^{\top})$ and using Lemma~\ref{lem:lipschitz} gives
\begin{equation}\label{eq:lipschitzproof}
\begin{aligned}
&\underset{\Q\in\B}{\sup}\Exp_{\Q}[(\ell(g,\xi)+\tau)_+]\\
&\leq \begin{cases}
\inf_{\lambda\geq 0}\lambda\eps+\frac{1}{N} \sum \limits_{i=1}^{N} (\ell(g,\hat{\xi}^{(i)})+\tau)_+ &\textup{if }L_{\textup{con}}\|g\|_{r,\ast}\leq \lambda \\
\infty &\textup{otherwise}
\end{cases}
\end{aligned}
\end{equation}
Hence, by resolving the $\inf$ and resorting to an epigraph formulation, we arrive at
\[
\begin{aligned}
&\inf_{\tau\in\real} -\tau\alpha+\sup_{\Q\in\B}\Exp_{\Q}[( \ell(g,\xi)+\tau)_+] \\
&\leq \begin{cases}
\inf_{\tau\in\real,s_i\in\real} -\tau\alpha+\eps L_{\textup{con}}\|g\|_{r,\ast} +\frac{1}{N}\sum_{i=1}^N s_i \\
\textup{s.t. } \tau+\ell(g,\hat{\xi}^{(i)}) \leq s_i\\
s_i\geq 0\\
\forall i\leq N
\end{cases}
\end{aligned}
\]
Substituting the definition for $\ell$ and recalling from~\eqref{eq:cvar} that we wish for the right hand side of the above inequality to be nonnegative gives
\[
\left\{ g\;\middle\vert
\begin{array}{l}
\inf_{\tau,s_i} -\tau \alpha + L_{\textup{con}}\eps\|g\|_{r,\ast} \\
\quad\quad\quad\quad\quad\quad\quad+\frac{1}{N}\sum_{i=1}^N s_i\leq 0\\
\tau + h(\Yfhat^{(i)} g) \leq s_i \\
s_i\geq 0 \\
\forall i\leq N
\end{array}
\right\} \subseteq G_{\textup{CVaR}} \,.
\]
We can remove the infimum from the above formulation by noting that there exist variables $\tau,s_i$ satisfying the constraints above if and only if the above infimum constraint holds.
The ``only if'' part is obvious. The ``if'' part can be split into two cases: the infimum is achieved, or the infimum is not achieved. If the infimum is achieved then the optimizer of the infimum satisfies the above. If the infimum is not achieved (i.e., the infimum is $-\infty$), then the first constraint is trivially satisfied, and one can find variables $\tau, s_i$ satisfying the remaining constraints. Note that~\eqref{eq:minmax} is an equality if $g\mapsto \ell(g,\xi)$ is convex for every $\xi$ and $\xi \mapsto \ell(g,\xi)$ is bounded on $\Xi$ for every $g$ (see,~\cite[Theorem 2.2]{AH-AC-JL:18}). Furthermore, ~\eqref{eq:lipschitzproof} is an equality if $\Xi=\real^{p(\Tini+\Tf)\lfloor \frac{T}{\Tini+\Tf}\rfloor}$. Since $h$ is convex, this proves the claim.
\end{IEEEproof}
We present an alternative result for the special case of piecewise affine constraints in which it is possible to obtain a tight set reformulation.
\begin{lemma}\longthmtitle{Piecewise Affine Constraint Reformulation}\label{lem:constraintreformulation_affine}
Let $\Xi=\{\xi\;\vert\; F\xi\leq d\}$ for some $F$ and $d$ of appropriate dimensions. Let $\ell(g,\xi)\defeq h(\xi_{p\Tini+i}^{\top} g,\dots,\xi_{p(\Tini+\Tf)}^{\top} g)$.
Assume $\ell$ is piecewise affine in $\xi$, i.e., $\ell(g,\xi)=\max_{k\leq K}\langle M_kg,\xi \rangle +b_k$ for matrices $M_k$, $b_k\in\real$ and $K\in\integerspositive$. Then 
\[
\left\{ g\;\middle\vert
\begin{array}{l}
\exists \tau,\lambda,s_i,\gamma_{ik}\; \textup{such that}\\
 -\tau \alpha + \lambda\eps +\frac{1}{N}\sum_{i=1}^N s_i\leq 0\\
 b_k+\tau + \langle M_kg,\hat{\xi}^{(i)}\rangle \\
 \qquad+\langle \gamma_{ik},d - F\hat{\xi}^{(i)} \rangle \leq s_i \\
\|F^{\top} \gamma_{ik} - M_k g\|_{r,\ast}\leq \lambda \\
s_i\geq 0 \\
\forall i\leq N,\; \forall k\leq K
\end{array}
\right\} \subseteq G_{\textup{CVaR}} \,.
\]
The sets coincide when $\ell(g,\xi)$ is bounded on $\Xi$ for all $g$.
\end{lemma}
\smallskip

\begin{IEEEproof}
Recall that the support set of the random variable $\xi$ is defined as $\Xi=\{\xi\;\vert\; F\xi\leq d\}$. Substituting $\ell(g,\xi)=\max_{k\leq K}\langle M_kg,\xi \rangle +b_k$ and applying~\cite[Corollary 5.1]{PME-DK:18} we have that
\[
\begin{aligned}
&\underset{\Q\in\B}{\sup}\Exp_{\Q}[(\ell(g,\xi)+\tau)_+]\\
&=\begin{cases}
\inf_{\lambda,s_i,\gamma_{ik}} \lambda\eps +\frac{1}{N}\sum_{i=1}^N s_i\\
\textup{s.t. } b_k+\tau + \langle M_kg,\hat{\xi}^{(i)}\rangle +\langle \gamma_{ik},d - F\hat{\xi}^{(i)} \rangle \leq s_i  \\
\|F^{\top} \gamma_{ik} - M_k g\|_{r,\ast}\leq \lambda \\
s_i\geq 0 \\
\forall i\leq N,\; \forall k\leq K
\end{cases}
\end{aligned}
\]
Invoking~\eqref{eq:cvar} and~\eqref{eq:minmax} and substituting the above yields the following set inclusion
\[
\left\{ g\; \middle\vert
\begin{array}{l}
\inf_{\tau,\lambda,s_i,\gamma_{ik}} -\tau \alpha + \lambda\eps +\frac{1}{N}\sum_{i=1}^N s_i\leq 0\\
\textup{s.t. } b_k+\tau + \langle M_kg,\hat{\xi}^{(i)}\rangle \\ 
\quad\qquad+\langle \gamma_{ik},d - F\hat{\xi}^{(i)} \rangle \leq s_i \\
\|F^{\top} \gamma_{ik} - M_k g\|_{r,\ast}\leq \lambda \\
s_i\geq 0 \\
\forall i\leq N,\; \forall k\leq K
\end{array}
\right\} \subseteq G_{\textup{CVaR}} \,. 
\]
Similarly as in the proof of Lemma~\ref{lem:constraintreformulation_lipschitz}, we can replace the infimum from the above formulation with existence of variables $\tau,\lambda,s_i, \gamma_{ik}$ satisfying the constraints above. Hence,
\[
\left\{ g\; \middle\vert
\begin{array}{l}
\exists \tau,\lambda,s_i,\gamma_{ik}\; \textup{such that}\\
 -\tau \alpha + \lambda\eps +\frac{1}{N}\sum_{i=1}^N s_i\leq 0\\
 b_k+\tau + \langle M_kg,\hat{\xi}^{(i)}\rangle \\
 \quad\quad+\langle \gamma_{ik},d - F\hat{\xi}^{(i)} \rangle \leq s_i  \\
\|F^{\top} \gamma_{ik} - M_k g\|_{r,\ast}\leq \lambda \\
s_i\geq 0 \\
\forall i\leq N,\; \forall k\leq K
\end{array}
\right\} \subseteq G_{\textup{CVaR}} \,.
\]
If $g\mapsto \ell(g,\xi)$ is convex for every $\xi$ and $\xi \mapsto \ell(g,\xi)$ is bounded on $\Xi$ for every $g$ then~\eqref{eq:minmax} reduces to equality (see,~\cite[Theorem 2.2]{AH-AC-JL:18}) which proves the claim.

\end{IEEEproof}
\begin{IEEEproof}[Proof of Theorem~\ref{thm:reformulation}]
Applying Lemma~\ref{lem:costreformulation} to the objective function and Lemma~\ref{lem:constraintreformulation_lipschitz} to the constraint of optimization problem~\eqref{eq:robustdeepc} gives the result.
\end{IEEEproof}
\begin{IEEEproof}[Proof of Theorem~\ref{thm:performanceguarantee}]
Choosing $\eps(\beta,N)$ as in~\eqref{eq:mineps} ensures that $\Prob^N(\Prob\in\B)\geq 1-\beta$. Hence, $J(g)$ is upper bounded by $\underset{\Q\in\B}{\sup} \quad \Exp_{\Q}\left[f_1(\Ufhat g) + f_2(\Yf g) + f_3(\Yp g -\yinihat) \right]$ with confidence $1-\beta$. Furthermore, by Theorem~\ref{thm:reformulation}, $J(g)$ is upper bounded by $\widehat{J}(g)$ with confidence $1-\beta$ for all $g$. Choosing $g=\hat{g}^{\star}$ proves the first statement.
From Lemma~\ref{lem:constraintreformulation_lipschitz}, the constraint set of~\eqref{eq:tractabledeepc} is a subset of the set of vectors $g$ satisfying $\underset{\Q\in\B}{\sup}\textup{CVaR}_{1-\alpha}^{\Q}(h(\Yf g))\leq 0$. Invoking the fact that $\Prob^N(\Prob\in\B)\geq 1-\beta$ proves the second statement.
\end{IEEEproof}
%

\bibliographystyle{IEEEtran}
\bibliography{IEEEabrv,JC}

\begin{IEEEbiography}[{\includegraphics[width=1in,height=1.25in,clip,keepaspectratio]{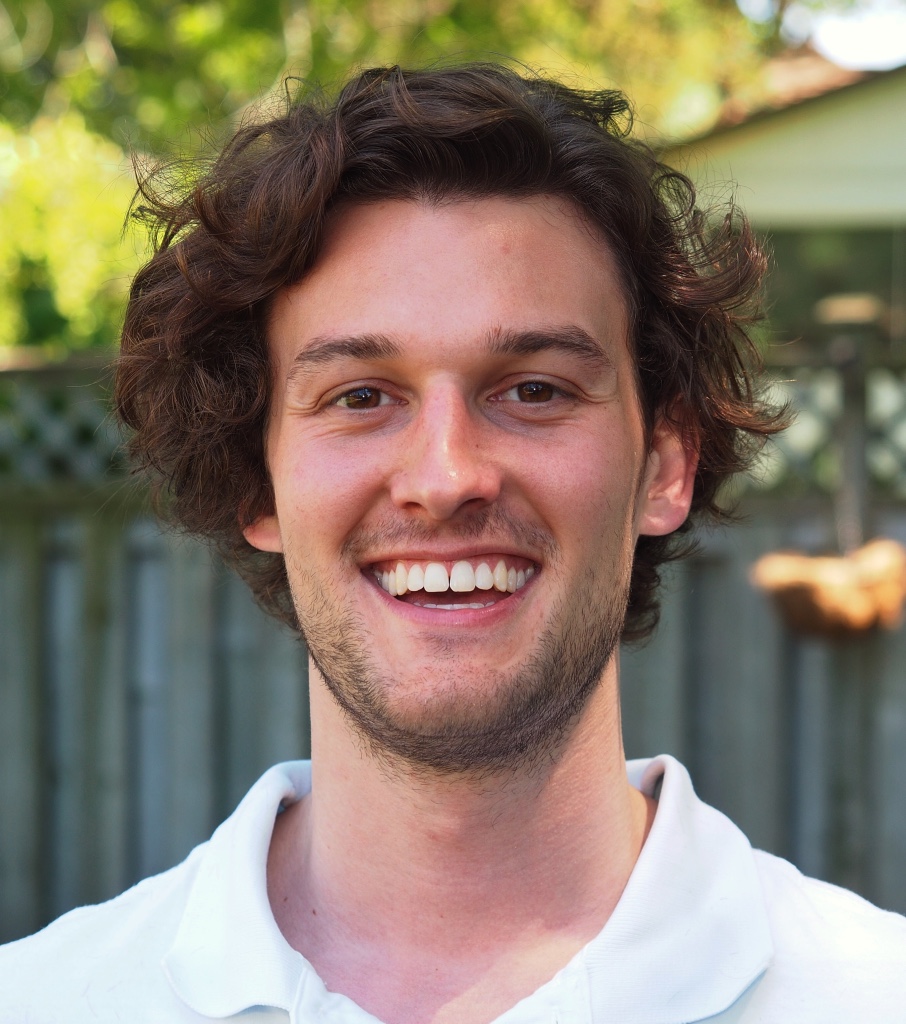}}]{Jeremy Coulson}
is a PhD student with the Automatic Control Laboratory at ETH Z\"{u}rich. He received his Master of Applied Science in Mathematics \& Engineering from Queen's University, Canada in August 2017. He received his B.Sc.Eng degree in Mechanical Engineering \& Applied Mathematics from Queen's University in 2015. His research interests include data-driven control methods, stochastic optimization, and control of partial differential equations.
\end{IEEEbiography}
%
\begin{IEEEbiography}[{\includegraphics[width=1in,height=1.25in,clip,keepaspectratio]{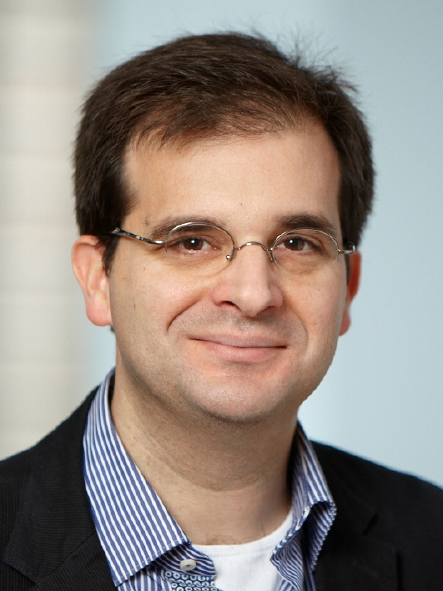}}]{John Lygeros}
completed a B.Eng. degree in electrical engineering in  1990 and an M.Sc. degree in Systems Control in 1991, both at Imperial College of Science Technology and Medicine, London, U.K.. In 1996 he obtained a Ph.D. degree from the Electrical Engineering and Computer Sciences Department, University of California, Berkeley. During the period 1996--2000 he held a series of post-doctoral researcher appointments at the Laboratory for Computer Science, M.I.T., and the Electrical Engineering and Computer Sciences Department at U.C. Berkeley. Between 2000 and 2003 he was a University Lecturer at the Department of Engineering, University of Cambridge, U.K., and a Fellow of Churchill College. Between 2003 and 2006 he was an Assistant Professor at the Department of Electrical and Computer Engineering, University of Patras, Greece. In July 2006 he joined the Automatic Control Laboratory at ETH Zurich, where he is currently serving as the Head of the Automatic Control Laboratory. His research interests include modelling, analysis, and control of hierarchical, hybrid, and stochastic systems, with applications to biochemical networks, automated highway systems, air traffic management, power grids and camera networks. John Lygeros is a Fellow of the IEEE, and a member of the IET and the Technical Chamber of Greece; since 2013 he serves as the Treasurer of the International Federation of Automatic Control.
\end{IEEEbiography}
%
\begin{IEEEbiography}[{\includegraphics[width=1in,height=1.25in,clip,keepaspectratio]{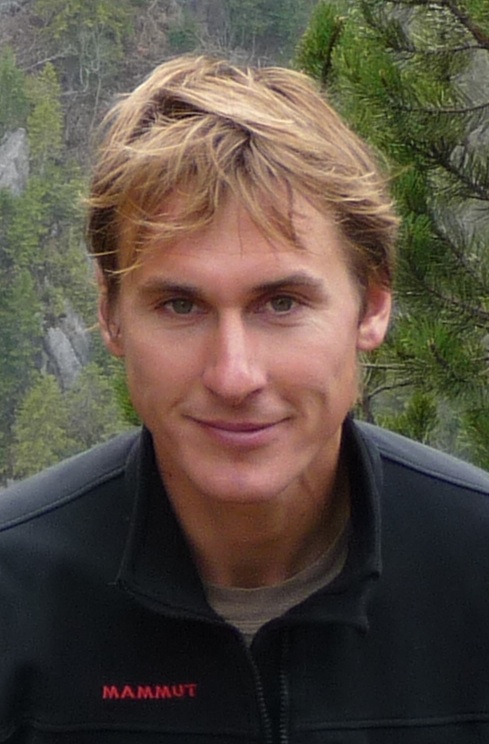}}]{Florian D\"orfler}
is an Associate Professor at the Automatic Control Laboratory at ETH Z\"urich. He received his Ph.D. degree in Mechanical Engineering from the University of California at Santa Barbara in 2013, and a Diplom degree in Engineering Cybernetics from the University of Stuttgart in 2008. From 2013 to 2014 he was an Assistant Professor at the University of California Los Angeles. His primary research interests are centered around control, optimization, and system theory with applications in network systems such as electric power grids, robotic coordination, and social networks. He is a recipient of the distinguished young research awards by IFAC (Manfred Thoma Medal 2020) and EUCA (European Control Award 2020). His students were winners or finalists for Best Student Paper awards at the European Control Conference (2013, 2019), the American Control Conference (2016), and the PES PowerTech Conference (2017). He is furthermore a recipient of the 2010 ACC Student Best Paper Award, the 2011 O. Hugo Schuck Best Paper Award, the 2012-2014 Automatica Best Paper Award, the 2016 IEEE Circuits and Systems Guillemin-Cauer Best Paper Award, and the 2015 UCSB ME Best PhD award.
\end{IEEEbiography}

\end{document}